\documentclass[11pt]{article}
\usepackage{graphicx,wrapfig}
\usepackage{graphicx}
\usepackage{flafter}
\usepackage{leftidx}
\usepackage{mathrsfs}
\usepackage{amssymb,amsmath}
\usepackage{bbding}
\usepackage{fancyhdr}
\usepackage{mathrsfs}
\usepackage{color}
\usepackage{stmaryrd}
\usepackage{amsfonts}
\usepackage{latexsym}
\usepackage{psfrag}
\usepackage{graphicx,subfigure}
\usepackage{fancyhdr,graphicx}
\usepackage{multicol}
\usepackage{dsfont}
\usepackage{bbm}
\usepackage{booktabs}
\usepackage{cite}
\usepackage{epstopdf}
\usepackage{diagbox}
\usepackage{booktabs}
\usepackage{multirow}
\usepackage{enumerate}
\usepackage{enumitem}
\allowdisplaybreaks

\date{}
\newcommand{\zd}{\,\mathrm{d}}

\newtheorem{theorem}{Theorem}[section]
\newtheorem{lemma}{Lemma}[section]
\newtheorem{remark}{Remark}
\newtheorem{example}{Example}[section]
\newtheorem{corollary}{Corollary}[section]
\numberwithin{equation}{section}

\newcommand{\diff}{\triangledown_{\!\tau}}
\newcommand{\abs}[1]{\left|#1\right|}

\newcommand{\abst}[1]{|#1|}
\newcommand{\bra}[1]{\left(#1\right)}
\newcommand{\brab}[1]{\big(#1\big)}
\newcommand{\braB}[1]{\Big(#1\Big)}
\newcommand{\brat}[1]{(#1)}

\newcommand{\kbrat}[1]{[#1]}
\newcommand{\myinner}[1]{\left\langle#1\right\rangle}
\newcommand{\myinnerb}[1]{\big\langle#1\big\rangle}
\newcommand{\myinnerB}[1]{\Big\langle#1\Big\rangle}
\newcommand{\mynorm}[1]{\left\|#1\right\|}
\newcommand{\mynormb}[1]{\big\|#1\big\|}

\newcommand{\mynormt}[1]{\|#1\|}

\newcommand{\hong}[1]{{\color{red}#1}}

\begin{document}
\title{Compatible $L^2$ norm convergence of variable-step L1 scheme 
for the time-fractional MBE mobel with slope selection}
\author{Yin Yang \thanks{
School of Mathematics and Computational Science, Xiangtan University,
Hunan National Applied Mathematics Center,
Xiangtan 411105, Hunan, China.
Yin Yang (yangyinxtu@xtu.edu.cn) is supported by
the National Natural Science Foundation of China Project 12071402,
the National Key Research and Development Program of China 2020YFA0713503,
the Project of Scientific Research Fund of the Hunan Provincial Science and
Technology Department 2020JJ2027.}
\quad Jindi Wang \thanks{
School of Mathematics and Computational Science, Xiangtan University,
Xiangtan 411105, Hunan, China.
Jindi Wang(wangjindixy@163.com) is supported by a grant
CX20200613 from Postgraduate Scientific Research Innovation Project of Hunan Province,
the Project of Scientific Research Fund of the Hunan Provincial Science and
Technology Department 2020ZYT003, 2018WK4006.}
\quad Yanping Chen\thanks{School of Mathematical Sciences, South China Normal University,
Guangzhou 510631, Guangdong, P.R. China; Yanping Chen (yanpingchen@scnu.edu.cn)
is supported by the State Key Program of National Natural Science Foundation
of China (11931003) and National Natural Science Foundation of China (41974133).}
\quad Hong-lin Liao\thanks{Corresponding author. ORCID 0000-0003-0777-6832. Department of Mathematics,
	Nanjing University of Aeronautics and Astronautics,
	Nanjing 211106, China; Key Laboratory of Mathematical Modelling
	and High Performance Computing of Air Vehicles (NUAA), MIIT, Nanjing 211106, China.
	Hong-lin Liao (liaohl@nuaa.edu.cn and liaohl@csrc.ac.cn)
	is supported by a grant 12071216 from
	National Natural Science Foundation of China.}
}
\maketitle
\normalsize

\begin{abstract}
The convergence of variable-step L1 scheme is studied for the time-fractional 
molecular beam epitaxy (MBE) model with slope selection.
A novel asymptotically compatible $L^2$ norm error estimate of the variable-step L1 scheme is established under a convergence-solvability-stability (CSS)-consistent time-step constraint. The CSS-consistent condition means that 
the maximum step-size limit required for convergence is of the same order to that for solvability and stability (in certain norms) as the small interface parameter $\epsilon\rightarrow 0^+$. To the best of our knowledge, it is the first time to establish such error estimate for nonlinear subdiffusion problems. The asymptotically compatible convergence means that the error estimate is compatible with that of backward Euler scheme for the classical MBE model as the fractional order $\alpha\rightarrow 1^-$.
Just as the backward Euler scheme can maintain the physical properties of the MBE equation, the variable-step L1 scheme can also preserve the corresponding properties of the time-fractional MBE model, including the volume conservation, variational energy dissipation law and $L^2$ norm boundedness. Numerical experiments are presented to support our  theoretical results.\\

\\
\noindent{\bf AMS subject classiffications.}\;\; 35Q99, 65M06, 65M12, 74A50
\end{abstract}
\section{Introduction}
\setcounter{equation}{0}

Consider the well-known Ehrlich--Schwoebel energy given as \cite{KohnYan:2003, MoldovanGolubovic:2000}
\begin{align}\label{eq:ESenergy}
E[\Phi]=\int_\Omega\frac{\epsilon^2}{2}\abst{\Delta \Phi}^2+F(\nabla \Phi)\zd{\mathbf{x}},
\end{align}
where the domain $\Omega=(0,L)^2\subset \mathbb{R}^2$,
the constant $\epsilon>0$ represents
the width of the rounded corners on the otherwise faceted crystalline thin films,
$\Phi$ is a scaled height function of a thin film, and
$F(\mathbf{v})=\frac{1}{4}(|\mathbf{v}|^2-1)^2$ is a nonlinear energy density function.
The MBE model with slope selection
can be viewed as the $L^2$ gradient flow
associated with the free energy \eqref{eq:ESenergy},
\begin{align}\label{eq:MBE gradient flow}
\partial_t \Phi :=-\kappa\mu\quad\text{with}\quad \mu:=\frac{\delta E}{\delta \Phi},
\end{align}
where 
$\mu$ is the vatiational derivative of the free energy $E$, 
$\kappa$ is a positive mobility constant, and 
the nonlinear vector functional $f(\mathbf{v})=F'(\mathbf{v})=(|\mathbf{v}|^2-1)\mathbf{v}$.
This model is widely used in material science because
it can accurately capture the growth of high-quality crystalline materials \cite{MoldovanGolubovic:2000}.
Under periodic boundary conditions, it is easy to check that the MBE system \eqref{eq:MBE gradient flow}
preserves the volume conservation  $\brab{\Phi(t),1}=\brab{\Phi(0),1}$,
the energy dissipation law
\begin{align}\label{eq:classical energy dissipation law}
\frac{\zd E}{\zd t}+\kappa\mynormb{\mu}_{L^2}^2=0,
\end{align}
and the following $L^2$ norm estimate, cf. the derivation of \eqref{eq:L2norm of continuous solution},
\begin{align}\label{eq:classical boundedness}
\mynorm{\Phi}_{L^2}^2\le\mynorm{\Phi_0}_{L^2}^2+
\frac{\kappa}{2}\abs{\Omega}t,
\end{align}
where $(\cdot, \cdot)$ and $\mynorm{\cdot}_{L^2}$ hong{are} the usual inner product and the associated $L^2$ norm.

Recently, many researchers paid great attention to
the time fractional phase field models \cite{ChenZhaoCaoWangZhang:2018, JiLiaoGongZhang:2018Adaptive, DuYangZhou:2020, MaskariKaraa:2021, TangYuZhou:2019, AinsworthMao: 2017}
to accurately describe the long time memory and the anomalously diffusive effects.
In this paper, we aim to develop a reliable numerical scheme for the time-fractional molecular beam epitaxy (TFMBE) model 
with slope selection, see 
\cite{ChenZhaoCaoWangZhang:2018, TangYuZhou:2019},
\begin{align}\label{eq:TFMBE}
\partial^\alpha_t \Phi=-\kappa\mu \quad\text{with}\quad
\mu=\epsilon^2\Delta^2 \Phi-\nabla\cdot f\brat{\nabla \Phi},
\end{align}
subject to the periodic boundary condition and initial condition
$\Phi(\mathbf{x}, 0):=\Phi_0(\mathbf{x})$.
As shown latter, this TFMBE model \eqref{eq:TFMBE} also retains some of continuous properties of the classical MBE model \eqref{eq:MBE gradient flow}. Here, $\partial^\alpha_t:={^C_0\mathcal{D}^\alpha_t}$ is the Caputo derivative of order $\alpha$,
\begin{align*}
{\partial^\alpha_t v}={^C_0\mathcal{D}^\alpha_t v}:=\mathcal{I}^{1-\alpha}_t v' \quad
\text{for $0<\alpha<1$,}
\end{align*}
where $\mathcal{I}^\beta_t$ is the fractional Riemann--Liouville integral operator of order $\beta>0$,
\begin{align*}
\brat{\mathcal{I}^\beta_t v}(t):=\int_0^t \omega_\beta(t-s)v(s)\zd{s} \quad
\text{with $\omega_\beta(t):=\frac{t^{\beta-1}}{\Gamma(\beta)}$}.
\end{align*}


\subsection{Continuous properties}
We describe some continuous properties of the TFMBE model \eqref{eq:TFMBE}, which are natural extensions of the physical properties of \eqref{eq:MBE gradient flow}, including the volume conservation, energy dissipation law \eqref{eq:classical energy dissipation law} and $L^2$ norm stability \eqref{eq:classical boundedness}.

Tang, Yu and Zhou \cite{TangYuZhou:2019} have established the volume conservation  $\brab{\Phi(t),1}=\brab{\Phi(0),1}$ and the following global energy dissipation law
\begin{align*}
	E\kbrat{\Phi(t)}\le E\kbrat{\Phi(0)} \quad \text{for $t>0$},
\end{align*}
which is quite different from the local energy decaying property
\eqref{eq:classical energy dissipation law}.
In order to be compatible with the classical model,
we consider a variational energy functional in \cite{LiaoTangZhou:2020Anenergy},
\begin{align}\label{eq:continuous variational energy definition}
	\mathcal{E}_\alpha[\Phi]:=E[\Phi]+\frac{\kappa}{2}\mathcal{I}_t^\alpha\mynormb{\mu}_{L^2}^2 \quad \text{for $t>0$}.
\end{align}
Obviously, this variational energy functional admits a local energy dissipation law
\begin{align}\label{eq:continuous variational energy dissipation law}
	\frac{\zd\mathcal{E}_\alpha}{\zd t}+\frac{\kappa}{2}\omega_\alpha(t)\mynormb{\mu}_{L^2}^2\le0.
\end{align}
This type energy functional $\mathcal{E}_\alpha[\Phi]$
is introduced firstly by Liao et al \cite{LiaoTangZhou:2020Anenergy}
in exploring the L1-type formula of Riemann--Liouville derivative
for the time-fractional Allen--Cahn equation.
If the fractional order $\alpha\rightarrow 1^-$,
the local energy decaying law \eqref{eq:continuous variational energy dissipation law}
asymptotically recovers the classical energy dissipation law in the form of
\begin{align*}
	\frac{\zd E}{\zd t}+\kappa\mynorm{\mu}_{L^2}^2\le0.
\end{align*}

In addition, by taking the $L^2$ inner product of the TFMBE model \eqref{eq:TFMBE} with $\Phi$, 
and using the Green's formula,
one gets
\begin{align}\label{eq:taking inner product of the TFMBE model}
	\bra{\partial^\alpha_t \Phi, \Phi}+\kappa\epsilon^2\mynorm{\Delta\Phi}_{L^2}^2-\kappa\bra{f\brat{\nabla \Phi}, \nabla\Phi}=0.
\end{align}
For the nonlinear term, one has
$$\bra{f\brat{\mathbf{v}}, \mathbf{v}}=
\bra{(|\mathbf{v}|^2-\tfrac1{2})^2-\tfrac14, 1}\ge -\frac14\bra{1,1}=-\frac14\abs{\Omega}.$$
By inserting it into \eqref{eq:taking inner product of the TFMBE model} 
and using the inequality
$\bra{\partial_t^\alpha \Phi, \Phi}
\ge \frac{1}{2}\partial_t^\alpha \mynorm{\Phi}_{L^2}^2$ from \cite[Lemma 2]{Alikhanov:2010Aprioriestimates}, 
one can reformulate the equation \eqref{eq:taking inner product of the TFMBE model}
into the following form
\begin{align*}
	\partial_t^\alpha \mynorm{\Phi}_{L^2}^2\le\frac{\kappa}{2}\abs{\Omega}, \quad t>0.
\end{align*}
By acting the Riemann-Liouville integral operator $\mathcal{I}^\alpha_t$ on both sides, one has
\begin{align}\label{eq:L2norm of continuous solution}
	\mynorm{\Phi}_{L^2}^2\le\mynorm{\Phi_0}_{L^2}^2+\frac{\kappa}{2}\abs{\Omega}\omega_{1+\alpha}(t),
	\quad t>0.
\end{align}
It is seen that, in the fractional order limit $\alpha\rightarrow 1^-$, the $L^2$ norm stability estimate
\eqref{eq:L2norm of continuous solution} is asymptotically compatible with 
 \eqref{eq:classical boundedness} of the  classical MBE equation \eqref{eq:MBE gradient flow}.

\subsection{Our contribution}

Some numerical methods were also proposed recently in \cite{ChenZhaoCaoWangZhang:2018,JiLiaoGongZhang:2018Adaptive,TangYuZhou:2019} for the TFMBE equation.
The numerical scheme in \cite{ChenZhaoCaoWangZhang:2018} utilized the fast L1 algorithm
for the Caputo derivative, 
but the $2-\alpha$ order of convergence was verified only experimentally. 
Ji \emph{et al.} \cite{JiLiaoGongZhang:2018Adaptive} suggested
a variable-step L$1^+$ scheme for the Caputo derivative
with second-order accuracy, and developed two Crank-Nicolson-type methods
based on the energy quadratization strategy. 
However, due to the lack of solution estimate, 
no convergence results are available in the literature for the numerical solutions of the TFMBE equation \eqref{eq:TFMBE}. 
In this paper, a rigorous $L^2$ norm convergence analysis
is presented for the variable-step L1 scheme. This scheme is asymptotically compatible with the backward Euler scheme for the classical MBE model \eqref{eq:MBE gradient flow} as the fractional order $\alpha\rightarrow 1^-$. Just as the backward Euler scheme can maintain the physical properties of the MBE equation, the variable-step L1 scheme can also preserve the corresponding properties of the time-fractional MBE model, including the volume conservation, varitional energy dissipation law
\eqref{eq:continuous variational energy dissipation law} and $L^2$ norm stability \eqref{eq:L2norm of continuous solution} at the discrete levels. 



\begin{table}[htbp]\label{table: time-step restriction}
	\renewcommand\arraystretch{2}
	\centering
	\caption{The CSS-consistent time-step conditions.}
	\vspace{0.3cm}
	\begin{tabular}{c|c|c|c}
		\hline
		\multicolumn{2}{c|}{~}&variable-step L1 scheme &backward Euler scheme ($\alpha\rightarrow1^-$)\\
		\hline
		\multicolumn{2}{c|}{Convergence}&$\tau_n\le\sqrt[\alpha]{2\omega_{2-\alpha}(1)\epsilon^2/\kappa}$&$\tau_n\le2\epsilon^2/\kappa$\\
		\hline
		\multicolumn{2}{c|}{Solvability}&$\tau_n\le\sqrt[\alpha]{4\omega_{2-\alpha}(1)\epsilon^2/\kappa}$&$\tau_n\le4\epsilon^2/\kappa$\\	
		\hline
		\multicolumn{2}{c|}{Energy stability}&$\tau_n\le\sqrt[\alpha]{4\omega_{2-\alpha}(1)\epsilon^2/\kappa}$&$\tau_n\le4\epsilon^2/\kappa$		\\
		\hline
		\multicolumn{2}{c|}{$L^2$ norm stability}&$\tau_n=O\bra{1}$&$\tau_n=O\bra{1}$\\
		\hline           
	\end{tabular}
	\vspace{0.3cm}
\end{table}

Many effective numerical methods \cite{Eyre:1998Unconditionally, ChenCondeWangWangWise:2012,FengWangiseZhang:2018,GongZhao:2019, JiLiaoZhang:2020Simple,	 ShenYang:2010,ShenXuYang:2018Thescalar, ShenXuYang:2019Anewclass,WangWangWise:2010Unconditionally,XuTang:2006}, including convex splitting methods, stabilized semi-implicit methods, exponential time differencing approaches and energy quadratization methods,    
have been explored rigorously for nonlinear phase field equations including the MBE model.
However, compared with the somewhat weak (or no) time-step constraints for solvability or the energy dissipation law,
the associated convergence analyses always suffer from very severe step-size restrictions with respect to the small interface parameter $\epsilon$ in the existing works. For example, the stablized method in \cite{FengWangiseZhang:2018} is unconditional energy stable with the step-size $\tau=O(1)$, but the convergence requires very small time-steps, nearly $\tau=O(\epsilon^{14})$.
It is an obvious defect at least in theoretical manner.
By making full use of the convexity of nonlinear functional $f(\mathbf{v})$, we establish an asymptotically compatible $L^2$ norm error estimate of the variable-step L1 scheme under a convergence-solvability-stability (CSS)-consistent time-step constraint. The CSS-consistent condition means that 
the maximum step-size limit required for convergence is of the same order to that for solvability and stability as the small interface parameter $\epsilon\rightarrow 0^+$. To the best of our knowledge, it is the first time to establish such error estimate for nonlinear subdiffusion problems.
Also, the imposed CSS-consistent time-step condition is asymptotically compatible with the time-step constraint of the backward Euler scheme as the fractional order $\alpha\rightarrow1^-$, see 
Table \ref{table: time-step restriction}.

In summary, our contribution is three-fold:
\begin{itemize}
	\item [$\bullet$] 
	By making use of the convexity of nonlinear bulk, a rigorous $L^2$ norm error estimate of the varaibel-step L1 scheme is established, maybe at the first time, under a CSS-consistent time-step condition. This estimate is robust and asymptotically compatible with that of the backward Euler scheme for the classical MBE model as $\alpha\rightarrow1^-$.
	
\item [$\bullet$] 
The variable step L1 scheme is proven to preserve the volume conservation, the variational energy dissipation law and $L^2$ norm stability
so that it is practically reliable in long-time simulations.
\item [$\bullet$] Several numerical examples are included to show the accuracy and effectiveness of the variable-step L1 scheme with an adaptive time-stepping strategy.
\end{itemize}

The rest of the paper is organized as follows.
Next section presents the nonuniform L1 implicit scheme and the unique solvability. The asymptotically compatible $L^2$ norm convergence is established in section 3. Section 4 addresses the discrete counterparts of the varitional energy dissipation law \eqref{eq:continuous variational energy dissipation law} and $L^2$ norm stability \eqref{eq:L2norm of continuous solution} at the discrete levels.
Some numerical examples are included in the last section.


\section{The variable-step L1 scheme and solvability}
\setcounter{equation}{0}

\subsection{Nonuniform L1 formula}

The TFMBE model \eqref{eq:TFMBE} has multi-scale behavior in a rough-smooth-rough pattern,
especially at an early stage of epitaxial growth on rough surfaces. It is practically useful to adopt some adaptive time-stepping strategy in the  coarsening dynamics approaching the steady state. It is desirable to investigate the time approximation on a general class of time meshes.

Consider
$0=t_0<\cdots<t_{k-1}<t_k<\cdots<t_N=T$ for a finite $T>0$.
Let the variable time-steps $\tau_k:=t_k-t_{k-1}$ for $1\le k\le N$.
We use the maximum step size $\tau:=\max_{1\le k\le N}\tau_k$,
and the adjoint time-step ratios $r_k:=\tau_k/\tau_{k-1}$ for $2\le k\le N$.
Given a grid function $\{v_k\}^N_{k=0}$,
let $\triangledown_{\tau}v^{k}=v^k-v^{k-1}$ and
$\partial_\tau v^k:=\triangledown_{\tau}v^{k}/\tau_k$
for $k\ge 1$.
The nonuniform L1 formula of Caputo derivative reads, see \cite{LiaoMcleanZhang:2019,LiaoYanZhang:2018Unconditional},
\begin{align}\label{eq: L1-discrete formula}
	\bra{\partial_\tau^\alpha v}^n
	:=\sum_{k=1}^n a_{n-k}^{(n)}\diff v^{k}\quad\text{for $n\ge1$,}
\end{align}
where the discrete coefficients $a_{n-k}^{(n)}$ are defined by
\begin{align}\label{eq:L1-Coefficient}
	a_{n-k}^{(n)}
	:=\frac{1}{\tau_{k}}\int_{t_{k-1}}^{t_{k}}\omega_{1-\alpha}(t_n-s)\zd{s}
	\quad\text{for $1\le{k}\le{n}$}.
\end{align}

\begin{figure}[htb!]
	\centering
	\includegraphics[width=3.5in,height=1.6in]{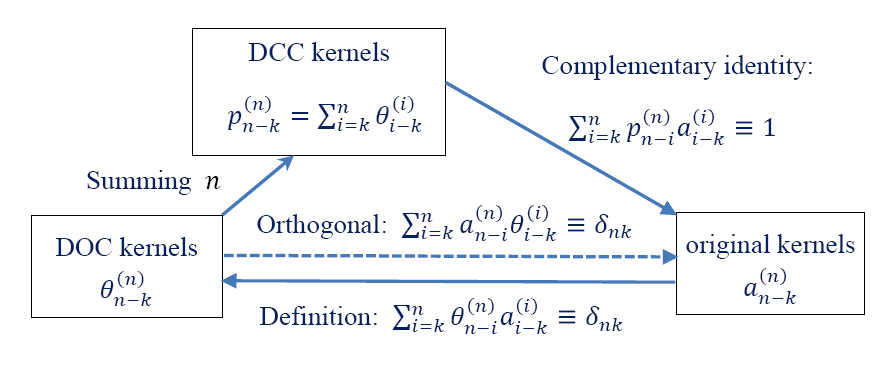}
	\caption{The relationship diagram of DOC and DCC kernels}
	\label{DOCDCCdiagram}
\end{figure}

We know that the discrete L1 kernels $a_{n-k}^{(n)}$ are positive and monotone on arbitrary time meshes \cite{LiaoTangZhou:2020Positive, LiaoZhuWang:2021}.
To deal with the discrete kernels, we introduce two important discrete tools, namely discrete orthogonal convolution (DOC) kernels and discrete complementary convolution (DCC) kernels.
The DOC kernels $\theta^{(n)}_{n-k}$ are defined via a recursive procedure \cite{LiaoZhang:2019Analysis}
\begin{align}\label{eq:DOC kernels definition}
	{\theta_{0}^{(n)}}:=\frac{1}{a_0^{(n)}}
	\quad \text{and} \quad
	{\theta_{n-k}^{(n)}}:=-\frac{1}{a_0^{(k)}}\sum^n_{j=k+1}{\theta_{n-j}^{(n)}}a^{(j)}_{j-k}
	\quad \text{for $1\le k\le n-1$}.
\end{align}
There has the following discrete orthogonal identity
\begin{align}\label{eq: discrete orthogonal identity}
	\sum_{j=k}^n \theta_{n-j}^{(n)} a_{j-k}^{(j)}\equiv\delta_{nk} \quad \text{for $1\le k\le n$},
\end{align}
where $\delta_{nk}$ is the Kronecker delta symbol.
The DCC kernels are defined as \cite{LiaoTangZhou:2020Positive}
\begin{align}\label{dcckernel}
	{p^{(n)}_{n-k}}:=\sum^n_{j=k}{\theta_{j-k}^{(j)}}
	\quad\text{for $1\le k\le n$}.
\end{align}
As proven in \cite[Subsection 2.2]{LiaoTangZhou:2020Positive}, the discrete convolution kernels ${p^{(n)}_{n-k}}$ are
complementary to the discrete L1 kernels $a_{n-k}^{(n)}$ in the following sense,
\begin{align}\label{eq: discrete complementary identity}
	\sum^n_{j=k}{p^{(n)}_{n-j}}a_{j-k}^{(j)}\equiv1
	\quad\text{for $1\le k\le n$}.
\end{align}
Figure \ref{DOCDCCdiagram} describes the above connections among three types of discrete convolution kernels. 

In the following convergence and stability analysis, we need the following result.
\begin{lemma}\cite[Lemma 2.1]{LiaoMcleanZhang:2019}\label{lemma:DCC kernels properties}
	For any $n\ge 2$, the DCC kernels $p_{n-k}^{(n)}$ in \eqref{dcckernel} satisfy,
	\begin{align*}
		p_{n-k}^{(n)}\ge 0 \quad\text{for $1\le k\le n$}, \quad \text{and}\quad
		\sum_{j=1}^n p_{n-j}^{(n)}\le \omega_{1+\alpha}(t_n).
	\end{align*}
\end{lemma}

\subsection{Fully discrete scheme}

Fourier pseudo-spectral method in space is adopted here.
Consider the discrete spatial grid
$\bar{\Omega}_h:=\{\mathbf{x}_h=(ih, jh)~|~0\le i,j\le M\}$ and $\Omega_h:=\bar{\Omega}_h\cap \Omega$,
where $M$ is an even positive integer and the uniform length $h:=L/M$.
Let $\mathscr{F}_M$ be the trigonometric polynomials space
(all trigonometric polynomials of degree up to $M/2$). Let $P_M: L^2(\Omega)\rightarrow \mathscr{F}_M$ and $I_M: L^2(\Omega)\rightarrow \mathscr{F}_M$
be the $L^2$-projection operator and the trigonometric interpolation operator of the periodic
function $v(\mathbf{x})\in L^2(\Omega)$, respectively, that is,
\begin{align*}
(P_M v)(\mathbf{x})=\sum_{m, n=-M/2}^{M/2-1}\hat{v}_{m, n}e_{m, n}(\mathbf{x}),\quad
(I_M v)(\mathbf{x})=\sum_{m, n=-M/2}^{M/2-1}\tilde{v}_{m,n}e_{m,n}(\mathbf{x}),
\end{align*}
where the basis function $e_{m,n}(\mathbf{x}):=e^{i\nu\bra{mx+ny}}$ with $\nu=2\pi/L$,
the coefficients $\hat{v}_{m, n}$ denote the standard Fourier coefficients of $v(\mathbf{x})$,
and the pseudo-spectral coefficients $\tilde{v}_{m,n}$ are determined
such that $(I_M v)(\mathbf{x}_h)=v(\mathbf{x}_h)$.
In turn, the Fourier pseudo-spectral
first- and second-order derivatives of $v_h$ are given by
\begin{align*}
\mathcal{D}_x v_h:=\sum_{m,n=-M/2}^{M/2-1} (i\nu m)\tilde{v}_{m,n} e_{m,n}(\mathbf{x}), \quad
\mathcal{D}_x^2 v_h:=\sum_{m,n=-M/2}^{M/2-1} (i\nu m)\tilde{v}_{m,n} e_{m,n}(\mathbf{x}).
\end{align*}
The notations $\mathcal{D}_y$ and $\mathcal{D}_y^2$ would be defined silimarly.
Accordingly, the discrete gradient $\nabla_h$ and Laplacian $\Delta_h$
in the point-wise sense are given by
\begin{align*}
\nabla_h v_h:=(\mathcal{D}_x v_h, \mathcal{D}_y v_h)^T\quad \text{and}
\quad \Delta_h v_h:=(\mathcal{D}^2_x v_h, \mathcal{D}^2_y v_h).
\end{align*}

In the numerical analysis, let
$\mathbb{V}_{h}:=\{v\,|\,v=\bra{v_h}\; \text{is L-periodic for}\; \mathbf{x}_h\in\bar{\Omega}_h\}$
be  the space of L-periodic grid functions.
For any functions $v, w\in \mathbb{V}_h$, the following discrete Green's formulas hold
$\myinner{-\Delta_h v, w}=\myinner{\nabla_h v, \nabla_h w}$ and
$\myinner{\Delta^2_h v, w}=\myinner{\Delta_h v, \Delta_h w}$. Also, we define the discrete inner product
$\myinner{v, w}:=h^2\sum_{\mathbf{x}_h\in\Omega_h} v_hw_h$,
the associated $L^2$ norm $\mynorm{v}:=\sqrt{\myinner{v, v}}$
and the discrete $L^p$ norm
$\mynorm{v}_{l^p}:=\sqrt[p]{h^2\sum_{\mathbf{x}_h\in\Omega_h}\abs{v_h}^p}$
for any grid functions $v, w \in \mathbb{V}_h$.
The discrete $H^1$ and $H^2$ norms are defined as
\begin{align*}
	\mynorm{v}^2_{H_h^1}:=\mynorm{v}^2+\mynorm{\nabla_h v}^2,\quad
	\mynorm{v}^2_{H_h^2}:=\mynorm{v}^2_{H_h^1}+\mynorm{\Delta_h v}^2.
\end{align*}

 We compute the numerical solution
$\phi_h^n\in \mathbb{V}_h$ of the TFMBE model \eqref{eq:TFMBE} by the fully implicit time-stepping scheme
\begin{align}\label{eq: nonuniform L1 scheme}
\bra{\partial_\tau^\alpha \phi_h}^n=-\kappa \mu_h^n
\quad\text{with}\quad
\mu_h^n=\epsilon^2\Delta_h^2 \phi_h^n-\nabla_h \cdot f\bra{\nabla_h \phi_h^n},
\end{align}
with the initial data $\phi_{h}^{0}=(P_M\Phi_{0})(\mathbf{x}_{h})$ for
$\mathbf{x}_{h}\in\Omega_{h}$. In order to facilitate our comparisons,
we also describe the backward Euler scheme for the calssical MBE model \eqref{eq:MBE gradient flow},
\begin{align}\label{eq:backward Euler for MBE}
\partial_\tau \phi_h^n=-\kappa \mu_h^n
\quad\text{with}\quad
\mu_h^n=\epsilon^2\Delta_h^2 \phi_h^n-\nabla_h \cdot f\bra{\nabla_h \phi_h^n}.
\end{align}
It is not difficult to check that, if the time-step size
$\tau_n\le {4\epsilon^2}/{\kappa}$,
the backward Euler scheme
\eqref{eq:backward Euler for MBE} is uniquely solvable  and fulfills the following energy dissipation law \cite{XuLiWu:2019}
\begin{align}\label{eq: BE energy dissipation law}
	\partial_\tau E[\phi^n]+\frac{\kappa}{2}\mynorm{\mu^n}^2 \le 0.
\end{align}

When the fractional index $\alpha\rightarrow 1^-$,
the discrete L1  kernels in \eqref{eq:L1-Coefficient} satisfy
$a_0^{(n)}\rightarrow 1/\tau_n$ and
$a_{n-k}^{(n)}\rightarrow 0$ for $1\le k\le n-1$.
Then, $\bra{\partial_\tau^\alpha \phi_h}^n\rightarrow \partial_\tau \phi_h^n$ as $\alpha\rightarrow 1^-$.
The nonuniform L1 scheme
\eqref{eq: nonuniform L1 scheme} is asymptotically compatible with the backward Euler scheme
\eqref{eq:backward Euler for MBE} in the fractional order limit $\alpha\rightarrow 1^-$.

\subsection{Unique solvability}

The full discrete scheme \eqref{eq: nonuniform L1 scheme} is volume conservative
and unique solvable.

\begin{lemma}\label{lemma:volume conservative}
	The full discrete scheme \eqref{eq: nonuniform L1 scheme} satisfies $\myinnerb{\phi^k, 1}=\myinnerb{\phi^0, 1}$ for $1\le k\le N$.
\end{lemma}
\begin{proof}The discrete Green's formula gives $\myinnerb{\mu^n, 1}=0$ from the second equation of \eqref{eq: nonuniform L1 scheme}. Thus the first equation of \eqref{eq: nonuniform L1 scheme}  yields $0=\myinnerb{\bra{\partial_\tau^\alpha \phi}^k,1}$ for $k\ge1$.
	Multiplying both sides of the above equality by the DOC kernels $\theta_{n-k}^{(n)}$ and
	summing $k$ from $k=1$ to $n$, we have
	\begin{align*}
		0=\myinnerB{\sum_{k=1}^n\theta_{n-k}^{(n)}\sum_{j=1}^ka_{k-j}^{(k)}\diff\phi^j,1}=
		\myinnerb{\diff\phi^n,1}\quad\text{for $n\ge1$,}
	\end{align*}
	where the summation order was exchanged and
	the discrete orthogonal identity \eqref{eq: discrete orthogonal identity} was applied
	in the second equality.
	It gives $\myinner{\phi^n, 1}=\myinner{\phi^{n-1}, 1}$
	and completes the proof.
\end{proof}


\begin{theorem}\label{theorem:unique solvable}
	If the time-step size satisfies
	\begin{align}\label{eq:unique solvable step size}
	\tau_n\le\sqrt[\alpha]{4\omega_{2-\alpha}(1)\epsilon^2/\kappa},
	\end{align}
	the nonuniform L1 scheme \eqref{eq: nonuniform L1 scheme} is uniquely solvable.
\end{theorem}
\begin{proof}
	We use the minimum principle of convex functional with a subspace of $\mathbb{V}_h$, that is,  $\mathbb{V}^*_h:=\big\{z\in \mathbb{V}_h\,|\,\myinnerb{z, 1}=\myinnerb{\phi^{n-1}, 1}\big\}$.
	Consider a discrete functional $G[z]$ on the space $\mathbb{V}^*_h$,
	\begin{align*}
	G[z]:=\frac{a_0^{(n)}}{2}\mynorm{z-\phi^{n-1}}^2
	+\myinnerb{\mathcal{L}^{n-1}, z-\phi^{n-1}}
	+\frac{\epsilon^2\kappa}{2}\mynorm{\Delta_h z}^2
	+\frac{\kappa}{4}\mynorm{\nabla_h z}_{l^4}^4-\frac{\kappa}{2}\mynorm{\nabla_h z}^2,
	\end{align*}
	where $n\ge 1$ and 
	$\mathcal{L}^{n-1}:=\sum_{k=1}^{n-1}a_{n-k}^{(n)}\diff  \phi^k.$
	This functional $G[z]$ is strictly convex under the time-step condition \eqref{eq:unique solvable step size} or $a_0^{(n)}\ge\kappa/\brat{4\epsilon^2}$.
	In fact,
	for any $\psi_h\in\mathbb{V}^*_h$,
	\begin{align*}
	\left.\frac{\zd^2 G}{\zd s^2}[z+s \psi]\right|_{s=0}
	&=a_0^{(n)}\mynorm{\psi}^2+\epsilon^2\kappa\mynorm{\Delta_h \psi}^2
	+3\kappa\mynorm{\nabla_h z\cdot\nabla_h \psi}^2-\kappa\mynorm{\nabla_h \psi}^2\\
	&\ge a_0^{(n)}\mynorm{\psi}^2+\epsilon^2\kappa\mynorm{\Delta_h \psi}^2
	+\kappa\myinnerb{\psi,\Delta_h\psi}\\
	&\ge \brab{a_0^{(n)}-\frac{\kappa}{4\epsilon^2}}\mynorm{\psi}^2\ge 0,
	\end{align*}
	where the Cauchy--Schwarz inequality and Young's inequality have been used in third step.
	Next, we show that the functional $G[z]$ is coercive on $\mathbb{V}^*_h$,
	\begin{align*}
	G[z]&\ge \frac{a_0^{(n)}}{2}\mynorm{z-\phi^{n-1}}^2+\myinner{\mathcal{L}^{n-1}, z-\phi^{n-1}}+\frac{\kappa}{4}\mynorm{\nabla_h z}_{l^4}^4
	-\frac{\kappa}{2}\mynorm{\nabla_h z}^2\\
	&\ge \kappa\mynorm{\nabla_h z}^2-\frac{1}{2a_0^{(n)}}\mynorm{\mathcal{L}^{n-1}}^2
	-\frac{9\kappa}{4}\abs{\Omega_h},
	\end{align*}
	where the inequality
	$\mynorm{\nabla_h z}_{l^4}^4\ge 6\mynorm{\nabla_h z}^2-9\abs{\Omega_h}$,
	due to the fact $\myinnerb{(\abst{\nabla_h z}^2-3)^2,1}\ge0$, was used in the last step.
 	Thus the functional $G[z]$ exists a unique minimizer, denote by $\phi_h^n$,
if and only if it solves the following equation	\begin{align*}
	\left.\frac{\zd G}{\zd s}[z+s \psi]\right|_{s=0}
	=\myinnerB{a_0^{(n)}(z-\phi^{n-1})+\sum_{k=1}^{n-1}a^{(n)}_{n-k}\diff  \phi^k
		+\kappa\epsilon^2\Delta_h z-\kappa\nabla_h\cdot f(\nabla_h z), \psi}=0.
	\end{align*}
	This equation holds for any $\psi_h\in\mathbb{V}^*_h$
	if and only if the unique minimizer $\phi_h^n\in\mathbb{V}^*_h$ solves
	\begin{align*}
	a_0^{(n)}\brat{\phi_h^n-\phi_h^{n-1}}+\sum_{k=1}^{n-1}a^{(n)}_{n-k}\diff  \phi_h^k
	+\kappa\epsilon^2\Delta_h \phi_h^n-\kappa\nabla_h\cdot f(\nabla_h \phi_h^n)=0,
	\end{align*}
	which is just the scheme \eqref{eq: nonuniform L1 scheme}.
	The proof is completed.
\end{proof}

	Note that, the time-step restriction \eqref{eq:unique solvable step size} of solvability 
	is sharp in the sense that it is compatible with that of the backward Euler scheme \eqref{eq:backward Euler for MBE}, that is, 
	\begin{align*}
	\tau_n\le\sqrt[\alpha]{4\omega_{2-\alpha}(1)\epsilon^2/\kappa}
	\quad\longrightarrow\quad \tau_n\le 4\epsilon^2/\kappa\qquad \text{as $\alpha\rightarrow 1^-$}.
	\end{align*}

\section{$L^2$ norm error estimate}
\setcounter{equation}{0}

This section presents the rigorous convergence analysis in the $L^2$ norm.
We use the standard semi-norms and norms of the Sobolev space $H^m(\Omega)$.
Let $\mathcal{C}_{per}^\infty(\Omega)$ be a set of
infinitely differentiable $L$-periodic functions defined on $\Omega$, and
$H_{per}^m(\Omega)$ be the closure of $\mathcal{C}_{per}^\infty(\Omega)$ in $H^m(\Omega)$,
endowed with the semi-norm $\abs{\cdot}_{H_{per}^m}$ and the norm $\mynorm{\cdot}_{H_{per}^m}.$
For the simplicity of notation, we denote $\abs{\cdot}_{H^m}:=\abs{\cdot}_{H_{per}^m}$,
$\mynorm{\cdot}_{H^m}:=\mynorm{\cdot}_{H_{per}^m}$,
and $\mynorm{\cdot}_{L^2}:=\mynorm{\cdot}_{H^0}$.

We recall the $L^2$-projection operator $P_M$ and
interpolation operator $I_M$
defined in Section 2,
and denote the $L^2$-projection
of exact solution $\Phi_M:=P_M \Phi$.
The following lemma lists the projection error $P_M v-v$, and the interpolation
error $I_M v-v$ in Sobolev space.
\begin{lemma}\cite{ShenTang:2006Spectral}\label{lemma:projectionerror}
	For any $v\in H^q_{per}(\Omega)$ and $0\le s\le q$, it holds that
	\begin{align*}
		\mynorm{P_M v-v}_{H^s}\le C h^{q-s}\abs{v}_{H^q},
		\quad \mynorm{P_M v}_{H^s}\le C\mynorm{v}_{H^s};
	\end{align*}
	and, in addition if $q>1,$
	\begin{align*}
		\mynorm{I_M v-v}_{H^s}\le C h^{q-s}\abs{v}_{H^q},
		\quad \mynorm{I_M v}_{H^s}\le C\mynorm{v}_{H^s}.
	\end{align*}
\end{lemma}

\subsection{Global consistency error}
Numerical tests in \cite{JiLiaoGongZhang:2018Adaptive} show that the TFMBE equation \eqref{eq:TFMBE} admits a weak singularity near the initial time,
like $\partial_t\Phi=O(t^{\alpha-1})$.
To complete the convergence analysis on nonuniform time meshes, it is reasonable to assume that,
\begin{align}\label{regularity}
\mynormb{\Phi}_{H^{m+4}}\le C_\phi,
\quad\mynormb{\partial_t^{\alpha}\Phi}_{H^{m}}\le C_\phi\quad\text{and}
\quad\mynormb{\partial_t^{(l)}\Phi}_{H^{m}}\le C_\phi(1+t^{\alpha-l}),
\end{align}
for $0<t\le T$ and $l=1, 2$,  where $m\ge 0$ is an integer,  $C_\phi$ denotes a generic positive constant.
Such a regularity assumption on the exact solution of time-fractional phase field models with the Caputo time derivative is standard in the numerical analysis\cite{DuYangZhou:2020, Stynes:2017, MaskariKaraa:2021, McLean:2020, JinLiZhou:2018}.

	The analytical solution of the TFMBE equation \eqref{eq:TFMBE}
	is weak singular at the initial time but regular away from the initial time.
We put a grading parameter $\gamma\ge 1$ and assume that
\begin{enumerate}[itemindent=1em]
	\item[\textbf{AG}.]
	there exists a constant $C_\gamma$, independent on the mesh,
	satisfies that the time-step sizes $\tau_k\le\tau \min\{1, C_\gamma t_k^{1-1/\gamma}\}$ for $1 \le k \le N$
	and $t_k \le C_\gamma t_{k-1}$ for $2 \le k \le N$.
\end{enumerate}
If the parameter $\gamma=1$, that means the mesh is quasi-uniform. As $\gamma$ increases, the initial step sizes are graded-like and become smaller compared to the others. On the other side, the assumption \textbf{AG} restricts only the maximum step size for the time mesh away from the initial time,
so that the step sizes can be adjusted according to the solution behaviors. This point is very important
in simulating the TFMBE model \eqref{eq:TFMBE} because it admits complex multi-scale behaviors in the long-time coarsening dynamics,
cf. Figures \ref{figure:etaEnergyEnalphaTimesteps} and \ref{figures:adaptive} in Section 5.

Let $\Upsilon^j=(\partial^\alpha_t v)(t_j)-(\partial_\tau^\alpha v)^j$
denote the local consistency error of the variable-step L1 formula \eqref{eq: L1-discrete formula}
at the time $t=t_j$.
We have the following results for the global convolution approximation error
$\sum_{j=1}^n p_{n-j}^{(n)}|\Upsilon^j|$, see
\cite[Lemma 3.1 and Lemma 3.3]{LiaoYanZhang:2018Unconditional}.

\begin{lemma}\cite[Lemma 3.1]{LiaoYanZhang:2018Unconditional}\label{lemma:consistency}
	For $v\in C^2(0,T]$ with $\int_0^T t|v_{tt}|\zd t<\infty$, 
		the global consistency error of the L1 formula \eqref{eq: L1-discrete formula} is bounded by
	\begin{align*}
	\sum_{j=1}^n p_{n-j}^{(n)}|\Upsilon^j|
	\le 2\sum_{k=1}^n p_{n-k}^{(n)}a_0^{(k)}\int^{t_k}_{t_{k-1}}(t-t_{k-1})|v_{tt}|\zd t.
	\end{align*}
\end{lemma}

We note that, the error bound in Lemma \ref{lemma:consistency} is valid on arbitrary time meshes and
is asymptotically compatible with the (global) truncation error of the backward Euler scheme
\eqref{eq:backward Euler for MBE}.  
Actually,  one has
\begin{align*}
\sum_{k=1}^n p_{n-k}^{(n)}a_0^{(k)}\int^{t_k}_{t_{k-1}}(t-t_{k-1})|v_{tt}|\zd t
\quad\longrightarrow\quad\sum_{k=1}^n\int^{t_k}_{t_{k-1}}(t-t_{k-1})|v_{tt}|\zd t
\quad\text{as $\alpha\rightarrow1^-$.}
\end{align*}
As desired, the limit is of temporal order $O(\tau)$.
On the other hand, the error bound in the next Lemma is not asymptotically compatible
in the fractional order limit $\alpha\rightarrow1^\hong{-}$.
This defect is mainly due to the lack of some proper estimates for the DCC kernels $p_{n-k}^{(n)}$;
however, it remains open to us up to now.

\begin{lemma}\cite[Lemma 3.3]{LiaoYanZhang:2018Unconditional}\label{lemma: gradeglobal}
	If  $v$ satisfies \eqref{regularity} and
	the meshes satisfy the assumption $\mathbf{AG}$,
	then the global consistency error of the L1 formula \eqref{eq: L1-discrete formula} can be bounded by
	\begin{align*}
		\sum_{j=1}^n p_{n-j}^{(n)}|\Upsilon^j|\le\frac{C_v}{\alpha(1-\alpha)}\tau^{\min\{2-\alpha, \gamma\alpha\}}~~\text{for}~~1\le n\le N.
	\end{align*}
\end{lemma}

\subsection{$L^2$ norm error estimate}


We are in a position to present the $L^2$ norm error estimate
for the variable-step L1 scheme \eqref{eq: nonuniform L1 scheme}.
The involving notation
$E_\alpha(z):=\sum_{k=0}^\infty\frac{z^k}{\Gamma(1+k\alpha)}$
denotes the Mittag--Leffler function.

\begin{theorem}\label{theorem:sharp L2 error estimate}
	Assume that the unique solution $\Phi$ of the TFMBE equation \eqref{eq:TFMBE}
	satisfies the regular condition \eqref{regularity}.
	If the time-step size 
	\begin{align}\label{eq:convergence step size}
	\tau_n\le \sqrt[\alpha]{2\omega_{2-\alpha}(1)\epsilon^2/\kappa},
	\end{align}
	then the numerical solution of the adaptive time-stpping  L1 scheme \eqref{eq: nonuniform L1 scheme}
	is unconditionally convergent in the discrete $L^2$ norm,
	\begin{align}\label{eq:convergence analysis}
	\mynorm{\Phi^n-\phi^n}
	&\le 2E_{\alpha}\braB{\frac{\kappa{t}_{n}^{\alpha}}{2\epsilon^2r_*}}
	\braB{C_\phi t_n^{\alpha} h^m+\max_{1\le j\le n}\sum_{k=1}^j p_{j-k}^{(j)}a_0^{(k)}\int^{t_k}_{t_{k-1}}\!\!(t-t_{k-1})\mynormb{\partial_{tt}\Phi}\zd t}
	\end{align}
	where $r_*:=\min_{1\le k\le N}\{1, r_k\}$ is the minimum step-ratio.
\end{theorem}

\begin{proof}
	We establish the error estimate for the fully implicit L1 scheme
	\eqref{eq: nonuniform L1 scheme}
	with the help of finite Fourier projection.
	The whole proof is divided into three steps.
	
	\textbf{Step1: Consistency error from projection (spatial discretization)}\quad
	Replacing the solution $\Phi$, the spatial operators $\Delta$ and $\nabla$
	with the projected solution $\Phi_M$, the discrete operators $\Delta_h$ and $\nabla_h$
	at the collocation points $\mathbf{x}_h\in\Omega_h$, respectively, one obtains
	\begin{align}\label{eq:projection continuous equation}
	\partial_t^\alpha \Phi_M\bra{\mathbf{x}_h, t}=-\kappa\epsilon^2\Delta_h^2 \Phi_M+\kappa\nabla_h\cdot f(\nabla_h \Phi_M)+\xi_h.
	\end{align}
Next, the $L^2$ norm of the consistency error $\xi_h$ will be evaluated.
		By subtracting \eqref{eq:TFMBE} from \eqref{eq:projection continuous equation}, 
		and applying the triangle inequality, one finds
	\begin{align}\label{eq:projection error bound}
	\mynorm{\xi}
	\le\mynorm{\partial_t^\alpha(\Phi-\Phi_M)}
	+\kappa\epsilon^2\mynorm{\Delta^2 \Phi-\Delta_h^2 \Phi_M}
	+\kappa\mynorm{\nabla\cdot f(\nabla \Phi)-\nabla_h\cdot f(\nabla_h \Phi_M)}.
	\end{align}
	Following the proof of \cite[Theorem 3.1]{LiaoJiZhang:2020pfc}, one can apply
	Lemma \ref{lemma:projectionerror} with the assumption \eqref{regularity} to find that
	\begin{align*}
	\mynorm{\Delta^2 \Phi-\Delta_h^2 \Phi_M}\le C_\phi h^m\quad\text{and}\quad
	\mynorm{\nabla\cdot f(\nabla \Phi)-\nabla_h\cdot f(\nabla_h \Phi_M)}\le C_\phi h^m.
	\end{align*}
	The projected time derivative $\partial_t^\alpha \Phi_M$ is the truncation of $\partial_t^\alpha \Phi$,
	for any $t>0$. Similarly, by using Lemma \ref{lemma:projectionerror} and the setting \eqref{regularity},
	one has
	\begin{align*}
	\mynorm{\partial^\alpha_t \brab{\Phi_M-\Phi}}\le C_{\phi}h^m\mynorm{\partial^\alpha_t\Phi}_{H^m}
	\le C_\phi h^m.
	\end{align*}
	In summary, we obtain that $\mynorm{\xi}\le C_\phi h^m$ for $t>0$ and then
	$$\mynorm{\xi(t_n)}\le C_\Phi h^m\quad \text{for $n\ge1$}.$$
		
	\textbf{Step2: Solution error from projection}\quad
	By replacing the numerical solution with the projection $\Phi_M^n(\mathbf{x}_h)$
	in the equation \eqref{eq: nonuniform L1 scheme}, one has
	\begin{align}\label{eq:projection numerical solution equation}
	(\partial^\alpha_\tau \Phi_M)^n
	=-\kappa\epsilon^2\Delta_h^2\Phi_M^n
	+\kappa\nabla_h\cdot f(\nabla_h \Phi_M^n)+\Upsilon_h^n+\xi_h^n
	\quad \text{for $n\ge 1$},
	\end{align}
	where $\Upsilon_h^n$ denotes the temporal consistency error, and
	$\xi_h^n:=\xi_h(t_n)$ is introduced from the projection equation \eqref{eq:projection continuous equation}.
	According to Lemma \ref{lemma:DCC kernels properties}, it is easy to derive that
	\begin{align}\label{eq-theorem:consistency from projection error}
	\sum_{j=1}^n p_{n-j}^{(n)}\mynorm{\xi^j}=
	\sum_{j=1}^n p_{n-j}^{(n)}\mynorm{\xi(t_j)}\le C_\phi \omega_{1+\alpha}(t_n)h^m
	\quad \text{for $n\ge 1$}.
	\end{align}
	Define $\Phi_M^n:= \Phi_M(\cdot, t_n)$.
	Let $e_h^n:=\Phi_M^n-\phi_h^n$ be the error between
	the finite Fourier projection $\Phi^n_M$ and the numerical solution $\phi^n$
	for any $\mathbf{x}_h\in\bar{\Omega}_h$.
	By subtracting the computational scheme \eqref{eq: nonuniform L1 scheme} from \eqref{eq:projection numerical solution equation},
	we get the following error system
	\begin{align*}
	(\partial^\alpha_\tau e_h)^n
	=-\kappa\epsilon^2\Delta_h^2e_h^n
	+\kappa\nabla_h\cdot \brab{f(\nabla_h \Phi_M^n)-f(\nabla_h \phi_h^n)}+\Upsilon_h^n+\xi_h^n,
	\end{align*}
	with the zero-valued data $e_h^0=0$.
	By taking the discrete inner product with $e^n$ and using the discrete Green's formula, one gets
	\begin{align}\label{eq-theorem:take inner product with error0}
	\myinner{(\partial^\alpha_\tau e)^n, e^n}
	=-\kappa\epsilon^2\mynorm{\Delta_he^n}^2
	+\kappa\mynorm{\nabla_h e^n}^2
	-\kappa\myinner{\mathbb{I}, \nabla_h e^n}+\myinner{\Upsilon^n+\xi^n, e^n},
	\end{align}
	where the nonlinear term
	$$\mathbb{I}:=|\nabla_h \Phi_M^n|^2\nabla_h \Phi_M^n-|\nabla_h \phi^n|^2\nabla_h \phi^n.$$ 
For any vectors $\mathbf{u}, \mathbf{v}\in \mathbb{R}^2$, it is not difficult to check that	
	\begin{align*}
	\myinnerb{|\mathbf{u}|^2\mathbf{u}-|\mathbf{v}|^2\mathbf{v}, \mathbf{u}-\mathbf{v}}=\frac1{2}\mynormb{|\mathbf{u}|^2-|\mathbf{v}|^2}^2+\frac1{2}\mynormb{\abs{\mathbf{u}-\mathbf{v}}(|\mathbf{u}|^2+|\mathbf{v}|^2)}^2\ge 0,
	\end{align*}
	which implies the nonlinear term $\myinner{\mathbb{I}, \nabla_he^n}\ge0$.
	Thus the equation \eqref{eq-theorem:take inner product with error0} reduces into
	\begin{align}\label{eq-theorem:take inner product with error}
		\myinner{(\partial^\alpha_\tau e)^n, e^n}
		\leq-\kappa\epsilon^2\mynorm{\Delta_he^n}^2
		+\kappa\mynorm{\nabla_h e^n}^2+
\myinner{\Upsilon^n+\xi^n, e^n},
	\end{align}
	For the term in left side of \eqref{eq-theorem:take inner product with error}, 
	by applying the decreasing property of the L1 kernels $a_{n-k}^{(n)}$,
	we get the following inequatity
	\begin{align}\label{eq-theorem:derivative inner product}
	\myinnerb{(\partial^\alpha_\tau e)^n, e^n}
	\ge\mynormb{e^{n}}\sum_{k=1}^{n}a_{n-k}^{(n)}\diff \mynormb{e^{k}}.
	\end{align}
	For the second term at the right side of \eqref{eq-theorem:take inner product with error},
	the Young's inequality also yields
	\begin{align}\label{eq-theorem:gradient inner product}
	\mynorm{\nabla_h e^n}^2\le
	\mynorm{\Delta_h e^n}\mynorm{e^n}
	\le \epsilon^2\mynorm{\Delta_h e^n}^2+\frac{1}{4\epsilon^2}\mynorm{e^n}^2.
	\end{align}
	 Inserting the above estimates \eqref{eq-theorem:derivative inner product}-\eqref{eq-theorem:gradient inner product} into \eqref{eq-theorem:take inner product with error}, we obtain
	\begin{align*}
	\mynormb{e^{n}}\sum_{k=1}^{n}a_{n-k}^{(n)}\diff \mynormb{e^{k}}\le \frac{\kappa}{4\epsilon^2}\mynorm{e^n}^{2}
	+\mynorm{\Upsilon^n}\mynorm{e^{n}}+\mynorm{\xi^n}\mynorm{e^{n}},
	\end{align*}
	which in turn gives the following estimate
	\begin{align*}
	\sum_{k=1}^{n}a_{n-k}^{(n)}\diff \mynormb{e^{k}}
	\le \frac{\kappa}{4\epsilon^2}\mynorm{e^n}+\mynorm{\Upsilon^n}+\mynorm{\xi^n}.
	\end{align*}
	Under the time-step restriction \eqref{eq:convergence step size},
	the well-known discrete fractional Gr\"{o}nwall inequality \cite[Theorem 3.2]{LiaoMcleanZhang:2019} yields
	\begin{align}\label{eq-theorem:error estimate}
	\mynorm{e^n}
	\le 2E_{\alpha}\braB{\frac{\kappa{t}_{n}^{\alpha}}{2\epsilon^2r_*}}
	\braB{C_\phi t_n^{\alpha} h^m+\max_{1\le j\le n}\sum_{k=1}^j p_{j-k}^{(j)}a_0^{(k)}\int^{t_k}_{t_{k-1}}(t-t_{k-1})\mynormb{\partial_{tt}\Phi}\zd t},
	\end{align}
where Lemma \ref{lemma:consistency}
and the bound \eqref{eq-theorem:consistency from projection error} were applied.
	
	\textbf{Step3: Error estimate}\quad
	Lemma \ref{lemma:projectionerror} gives
	the error of finite Fourier projection,
	\begin{align}\label{eq-theorem:the error of finite Fourier projection}
	\mynorm{\Phi^n_M-\Phi^n}=&\, \mynorm{I_M(\Phi^n_M-\Phi^n)}_{L^2}
	\le C_\phi \mynorm{\Phi^n-\Phi^n_M}_{L^2} \le C_\phi h^{m}\abs{\Phi^n}_{H^m}.
	\end{align}
	The triangle inequality with the estimates \eqref{eq-theorem:error estimate} and \eqref{eq-theorem:the error of finite Fourier projection} gives the claimed result.
\end{proof}

	The $L^2$ norm error estimate \eqref{eq:convergence analysis} is asymptotically compatible
	with  that of the backward Euler scheme \eqref{eq:backward Euler for MBE} in the limit $\alpha\rightarrow1^-$. As remarked for Lemma \ref{lemma:consistency}, we see that the error estimate \eqref{eq:convergence analysis} of the variable-step L1 scheme \eqref{eq: nonuniform L1 scheme} is $\alpha$-robust (not necessarily at the optimal convergence rate) in the sense of \cite{ChenStynes:2021}, in which an $\alpha$-robust bound was derived for the L1 formula. 
	Interested readers can follow the approach of \cite{ChenStynes:2021} to obtain the $\alpha$-robust estimate with optimal convergence order on graded meshes.
	 We emphasize that the presented $\alpha$-robust error estimate \eqref{eq:convergence analysis} is also 
	 mesh-robust for any finite $r_*$.

\begin{corollary}
	Assume that the unique solution $\Phi$ of the TFMBE equation \eqref{eq:TFMBE}
	satisfies the regular condition \eqref{regularity}.
	If the meshes satisfy $\mathbf{AG}$ and  \eqref{eq:convergence step size},  it holds that
	\begin{align*}
	\mynorm{\Phi^n-\phi^n}
	\le \frac{C_\phi}{\alpha(1-\alpha)}
	E_{\alpha}\braB{\frac{\kappa{t}_{n}^{\alpha}}{2\epsilon^2r_*}}
	\brab{t_n^{\alpha} h^m+\tau^{\min\{2-\alpha,\gamma\alpha\}}}\quad\text{for $1\leq{n}\leq{N}$}.
	\end{align*}
	The optimal accuracy is $O(\tau^{2-\alpha})$ if the grading parameter
	$\gamma\ge \max\{1, (2-\alpha)/\alpha\}$.
\end{corollary}

\section{Energy dissipation law and $L^2$ norm stability}
\setcounter{equation}{0}

The following lemma shows a discrete gradient structure
of the  L1 fromula \eqref{eq: L1-discrete formula}, which plays an important role
in the construction of discrete variational energy law.

\begin{lemma}\label{lemma:inequality about the L1 discrete convolution kernels and DCC kernels}
	For any real sequence $\{v_k\}_{k=1}^n$, it holds that
	\begin{align*}
		2v_n\sum^n_{j=1} {a^{(n)}_{n-j}}v_j
		\ge a_0^{(n)}v_n^2
		+\sum^n_{k=1} {p^{(n)}_{n-k}}\braB{\sum^k_{j=1} {a^{(k)}_{k-j}}v_j}^2
		-\sum^{n-1}_{k=1} {p^{(n-1)}_{n-1-k}}\braB{\sum^k_{j=1} {a^{(k)}_{k-j}}v_j}^2.
	\end{align*}
\end{lemma}

\begin{proof}From \cite[Lemma 2.4]{LiaoZhuWang:2021}, for any real sequence $\{w_k\}_{k=1}^n$, it holds
		\begin{align}\label{eq-lemma:inequality-2}
		2w_n\sum^n_{k=1} {\theta^{(n)}_{n-k}}w_k
		\ge
		\sum^n_{k=1} {p^{(n)}_{n-k}}w_k^2
		-\sum^{n-1}_{k=1} {p^{(n-1)}_{n-1-k}}w_k^2
		+\frac{1}{\theta_0^{(n)}}\braB{\sum^n_{k=1} {\theta^{(n)}_{n-k}}w_k}^2,
	\end{align}
where $\theta^{(n)}_{n-k}$ are the DOC kernels with respect to the L1 kernels $a^{(n)}_{n-j}$.
	We define
	$$v_j:=\sum_{k=1}^j \theta_{j-k}^{(j)} w_k.$$
	Multiplying both sides of this identity by the L1 kernels $a^{(n)}_{n-j}$ and
	summing $j$ from $j=1$ to $n$, we obtain
	\begin{align*}
		\sum_{j=1}^na_{n-j}^{(n)}v_j=\sum_{j=1}^na_{n-j}^{(n)}\sum_{k=1}^j \theta_{j-k}^{(j)} w_k
		=\sum_{k=1}^nw_k\sum_{j=k}^na_{n-j}^{(n)}\theta_{j-k}^{(j)}=w_n.
	\end{align*}
	The desired inequality is verified by inserting the above formulas of $v_j$ and $w_n$ into \eqref{eq-lemma:inequality-2}.
\end{proof}

We define a discrete counterpart of the variational energy
\eqref{eq:continuous variational energy definition} as follows
\begin{align}\label{eq:discrete variational energy definition}
\mathcal{E}_\alpha[\phi^0]:=E[\phi^0]\quad\text{and}\quad
\mathcal{E}_\alpha[\phi^n]:=E[\phi^n]+\frac{\kappa}{2}\sum^n_{j=1}{p^{(n)}_{n-j}}\mynormb{\mu^j}^2
\quad\text{for $n\ge 1$},
\end{align}
where $E[\phi^n]$ denotes the discrete counterpart of free energy \eqref{eq:ESenergy},
\begin{align*}
E[\phi^n]:=\frac{\epsilon^2}{2}\mynorm{\Delta_h \phi^n}^2
+\frac{1}{4}\mynormb{\abs{\nabla_h \phi^n}^2-1}^2.
\end{align*}
Here, the DCC kernels $p_{n-j}^{(n)}$ would be regarded as the discrete kernels
of the Riemann-Liouville fractional integral $\mathcal{I}_t^\alpha$, see \cite{LiaoMcleanZhang:2019},
$(\mathcal{I}_t^\alpha v)(t_n)\approx\sum_{j=1}^n p_{n-j}^{(n)} v^j$.
\begin{theorem}\label{Thm:Energystability}
Under the time step restriction \eqref{eq:unique solvable step size},
the L1 scheme \eqref{eq: nonuniform L1 scheme}
preserves the variational energy dissipation law at each time level,
\begin{align}\label{discrete energy dissipation law}
\partial_\tau\mathcal{E}_{\alpha}[\phi^n]\le 0\quad \text{for $1\le n\le N$}.
\end{align}
\end{theorem}
\begin{proof}
By taking the inner product of \eqref{eq: nonuniform L1 scheme}
with $\diff\phi^n$, it is easy to find
\begin{align}\label{eq-theorem:discrete variational energy dissipation law-1}
\myinnerB{\sum_{k=1}^n{a_{n-k}^{(n)}}\diff  \phi^k,\diff  \phi^n}
+\kappa\epsilon^2\myinner{\Delta_h \phi^n,\Delta_h\diff  \phi^n}
+\kappa\myinner{f(\nabla_h \phi^n),\nabla_h\diff  \phi^n}
=0.
\end{align}
For the first term on the left hand side,
by taking $v_k=\diff  \phi^k$ in Lemma \ref{lemma:inequality about the L1 discrete convolution kernels and DCC kernels}, we have
\begin{align*}
\myinnerB{\sum_{k=1}^n{a_{n-k}^{(n)}}\diff  \phi^k,\diff  \phi^n}
\ge\frac{\kappa^2}{2}\sum^n_{k=1} {p^{(n)}_{n-k}}\mynormb{\mu^k}^2
-\frac{\kappa^2}{2}\sum^{n-1}_{k=1} {p^{(n-1)}_{n-1-k}}\mynormb{\mu^k}^2
+\frac{a_0^{(n)}}{2}\mynorm{\diff  \phi^n}^2.
\end{align*}
By using Young's inequality, one has
$$\myinnerb{\abst{\nabla_h \phi^n}^2\nabla_h \phi^n,\nabla_h \phi^{n-1}}
\le \frac{3}{4}\mynorm{\nabla_h \phi^n}^4_{l^4}
+\frac{1}{4}\mynorm{\nabla_h \phi^{n-1}}^4_{l^4}.$$
Then the nonlinear term can be bounded by
\begin{align*}
\myinnerb{\abst{\nabla_h \phi^n}^2\nabla_h \phi^n,\nabla_h\diff  \phi^n}
&=\mynorm{\nabla_h \phi^n}^4_{l^4}-\myinnerb{\abst{\nabla_h \phi^n}^2\nabla_h \phi^n,\nabla_h \phi^{n-1}}\nonumber\\
&\ge\frac{1}{4}\brab{\mynorm{\nabla_h \phi^n}^4_{l^4}-\mynorm{\nabla_h \phi^{n-1}}^4_{l^4}}.
\end{align*}
Furthermore, the identity $2a(a-b)=a^2-b^2+(a-b)^2$ yields
\begin{align*}
\epsilon^2\myinnerb{\Delta_h \phi^n,\Delta_h\diff  \phi^n}
&=\frac{\epsilon^2}{2}\brab{\mynorm{\Delta_h \phi^n}^2-\mynorm{\Delta_h \phi^{n-1}}^2
+\mynorm{\Delta_h \diff  \phi^n}^2},\\
-\myinnerb{\nabla_h \phi^n,\nabla_h\diff  \phi^n}
&=\frac{1}{2}\brab{\mynorm{\nabla_h \phi^{n-1}}^2-\mynorm{\nabla_h \phi^n}^2
-\mynorm{\nabla_h\diff  \phi^n}^2}.
\end{align*}
Thus collecting the above estimates, it follows from  \eqref{eq-theorem:discrete variational energy dissipation law-1} that
\begin{align}\label{eq-theorem:discrete variational energy dissipation law-2}
\mathcal{E}_{\alpha}[\phi^n]-\mathcal{E}_{\alpha}[\phi^{n-1}]
+\frac{\epsilon^2}{2}\mynorm{\Delta_h \diff  \phi^n}^2
-\frac{1}{2}\mynorm{\nabla_h\diff  \phi^n}^2
+\frac{a_0^{(n)}}{2\kappa}\mynorm{\diff  \phi^n}^2
\le 0.
\end{align}
By using the Young's inequality, one gets
\begin{align*}
\mynorm{\nabla_h\diff  \phi^n}^2
\le\mynorm{\Delta_h\diff  \phi^n}\cdot\mynorm{\diff  \phi^n}
\le\epsilon^2
\mynorm{\Delta_h\diff  \phi^n}^2
+\frac{1}{4\epsilon^2}\mynorm{\diff  \phi^n}^2.
\end{align*}
Then we have
\begin{align}\label{eq-theorem:discrete variational energy dissipation law-3}
\mathcal{E}_{\alpha}[\phi^n]-\mathcal{E}_{\alpha}[\phi^{n-1}]
+\frac1{2\kappa}\brab{a_0^{(n)}-\frac{\kappa}{4\epsilon^2}}
\mynorm{\diff  \phi^n}^2\le 0.
\end{align}
Under the time-step restriction \eqref{eq:unique solvable step size}, the claimed inequality follows immediately.
\end{proof}

	Note that the DCC kernels satisfy
	$p_{n-j}^{(n)}\rightarrow \tau_j$ for $1\le j\le n$ as $\alpha\rightarrow 1^-$. Then one has
	\begin{align}
	\mathcal{E}_\alpha[\phi^n]\quad\longrightarrow\quad
	E[\phi^n]+\frac{\kappa}{2}\sum^n_{j=1}\tau_j\mynormb{\mu^j}^2
	\qquad \text{as $\alpha\rightarrow 1^-$}.
	\end{align}
	We see that the discrete variational energy dissipation law \eqref{discrete energy dissipation law}
	is asymptotically compatible with the classical energy law \eqref{eq: BE energy dissipation law}
	of the backward Euler scheme, that is,
	\begin{align*}
	\partial_\tau\mathcal{E}_{\alpha}[\phi^n]\le 0\quad\longrightarrow\quad
	\partial_\tau E[\phi^n]+\frac{\kappa}{2}\mynorm{\mu^n}^2 \le 0\qquad \text{as $\alpha\rightarrow 1^-$}.
	\end{align*}

\begin{theorem}\label{theorem: $L^2$ stability of numerical solution}
	The discrete solution $\phi^n$ of the variable-step L1 scheme 
	\eqref{eq: nonuniform L1 scheme}
	is unconditionally $L^2$ norm stable in the sense that
	\begin{align}\label{eq: boundedness of TFMBE}
	\mynormb{\phi^n}^2\le \mynorm{\phi^0}^2+\frac{\kappa}{2}\abs{\Omega_h}\omega_{1+\alpha}(t_n).
	\end{align}
\end{theorem}
\begin{proof}
	By taking the $L^2$ inner product of the nonuniform L1 scheme \eqref{eq: nonuniform L1 scheme} 
	with $\phi^n$, then adding up two results and 
	using the discrete Green's formula, we obtain
	\begin{align}\label{eq:boundedness-take inner product}
	\myinnerB{\sum_{k=1}^n a_{n-k}^{(n)}\diff \phi^k, \phi^n}
	+\kappa\epsilon^2\mynorm{\Delta_h \phi^n}^2+\kappa\myinner{f(\nabla_h \phi^n), \nabla_h\phi^n}=0.
	\end{align}
	One applies the decreasing property of  $a_{n-k}^{(n)}$ to get
	\begin{align}\label{eq:boundedness-Discrete derivative}
	\myinnerB{\sum_{k=1}^n a_{n-k}^{(n)}\diff \phi^k, \phi^n}
	\ge \frac{1}{2}\sum_{k=1}^n a_{n-k}^{(n)}\diff \mynormb{\phi^k}^2.
	\end{align}
	For the nonlinear term at the left hand side, it holds that
	\begin{align}\label{eq:boundedness-nonlinear term}
	\myinner{f(\nabla_h \phi^n), \nabla_h\phi^n}
	=\mynormb{\abst{\nabla_h\phi^n}^2-\tfrac{1}{2}}^2-\frac{1}{4}\abs{\Omega_h}
	\ge -\frac{1}{4}\abs{\Omega_h}.
	\end{align}
	Inserting above estimates \eqref{eq:boundedness-Discrete derivative} and \eqref{eq:boundedness-nonlinear term} into \eqref{eq:boundedness-take inner product}, one yields,
	\begin{align}\label{eq:boundedness-synthesize estimates}
	\sum_{k=1}^n a_{n-k}^{(n)}\diff \mynormb{\phi^k}^2
	-\frac{\kappa}{2}\abs{\Omega_h}\le 0.
	\end{align}
	We replace the index $n$ with $j$ in above inequality, 
	then multiply by $p_{n-j}^{(n)}$ and sum over $j$ from $1$ to $n$ to obtain
	\begin{align}\label{eq:boundedness-acting p and summing up}
	\sum_{j=1}^{n}p_{n-j}^{(n)}\sum_{k=1}^j a_{j-k}^{(j)}\diff \mynormb{\phi^k}^2
	- \frac{\kappa}{2}\abs{\Omega_h}\sum_{j=1}^{n}p_{n-j}^{(n)}\le 0.
	\end{align}
	By exchanging the order of summation, one applies
	the complementary identity \eqref{eq: discrete complementary identity} to get 
	\begin{align*}
	\sum_{j=1}^{n}p_{n-j}^{(n)}\sum_{k=1}^j a_{j-k}^{(j)}\diff \mynormb{\phi^k}^2
	=\sum_{k=1}^n\diff \mynormb{\phi^k}^2\sum_{j=k}^{n}p_{n-j}^{(n)}a_{j-k}^{(j)}
	=\sum_{k=1}^n\diff \mynormb{\phi^k}^2
	=\mynormb{\phi^n}^2-\mynormb{\phi^0}^2.
	\end{align*}
Thus, by using Lemma \ref{lemma:DCC kernels properties}, it follows from \eqref{eq:boundedness-synthesize estimates} that
	\begin{align*}
	\mynormb{\phi^n}^2\le\mynormb{\phi^0}^2
	+\frac{\kappa}{2}\abs{\Omega_h}\omega_{1+\alpha}(t_n).
	\end{align*}
	The proof is completed.
\end{proof}

	As the fractional order $\alpha\rightarrow1^{-}$, the $L^2$ norm boundedness \eqref{eq: boundedness of TFMBE} is asymptotically compatible
	with the $L^2$ norm solution estimate of backward Euler scheme, that is, 
	\begin{align*}
		\mynormb{\phi^n}^2\le \mynorm{\phi^0}^2+\frac{\kappa}{2}\abs{\Omega_h}t_n.
	\end{align*}
	This estimate can be derived by following the proof of Theorem \ref{theorem: $L^2$ stability of numerical solution}.
	
\begin{remark}
		Consider the convex splitting scheme \cite{Eyre:1998Unconditionally,JiLiaoZhang:2020Simple} for the TFMBE  model \eqref{eq:TFMBE},
		\begin{align*}
		\bra{\partial_\tau^\alpha \phi_h}^n=-\kappa \mu_h^n
		\quad\text{with}\quad
		\mu_h^n=\epsilon^2\Delta_h^2 \phi_h^n-\nabla_h \cdot \bra{|\nabla_h \phi_h^n|\nabla_h \phi_h^n}+\Delta_h \phi_h^{n-1}.
		\end{align*}
	    It is not difficult to check that this scheme is volume conservative and unconditionally solvable.  
	    With slight modifications to the proofs of Theorems \ref{Thm:Energystability} and \ref{theorem: $L^2$ stability of numerical solution}, one can show that the convex splitting scheme is unconditionally stable with respect to the discrete energy and the $L^2$ norm. 
	    That is to say, the time-step requirements for the solvability and stability are about $\tau_n=O(1)$. Nonetheless, the $\alpha$-robust, first-order convergence still requires the time-step restriction \eqref{eq:convergence step size}. In this case, the condition \eqref{eq:convergence step size} is not a
	    CSS-consistent time-step constraint.
	\end{remark}
	
\section{Numerical experiments}
\setcounter{equation}{0}

In this section,
we present several numerical examples to test
the accuracy and efficiency of the  L1 scheme
\eqref{eq: nonuniform L1 scheme} for the TFMBE model \eqref{eq:TFMBE}.
We use a simple fixed-point iteration algorithm with the termination error $10^{-12}$
to solve the resulting nonlinear equations at each time step.
Also, the sum-of-exponentials technique \cite{LiaoYanZhang:2018Unconditional} with the absolute tolerance error $\varepsilon=10^{-12}$
is employed to speed up the convolution computation of the L1 formula \eqref{eq: L1-discrete formula}.

\subsection{Convergence test}
We present an accuracy check for the L1 scheme
\eqref{eq: nonuniform L1 scheme}.
The time accuracy is focused on
and the spatial error
(standard spectral accuracy produced by the Fourier pseudo-spectral method) is negligible.
The experimental convergence order in time is computed by
\begin{align*}
\text{Order}:=\frac{\log\bra{e(N)/e(2N)}}{\log\bra{\tau(N)/\tau(2N)}},
\end{align*}
where the discrete $L^2$ norm error $e(N):=\max_{1\le n\le N}\mynorm{\Phi^n-\phi^n}$
and $\tau(N)$ denotes the maximum time-step size for total $N$ subintervals.

\begin{table}[htb!]
\begin{center}
\caption{Temporal error of \eqref{eq: nonuniform L1 scheme} for $\alpha=0.8$ with $\gamma_{\mathrm{opt}}=1.5$.}\label{table:StabilizedL1-InitialSingularity-1} \vspace*{0.3pt}
\def\temptablewidth{1.0\textwidth}
{\rule{\temptablewidth}{0.5pt}}
\begin{tabular*}{\temptablewidth}{@{\extracolsep{\fill}}cccccccccc}
\multirow{2}{*}{$N$} &\multirow{2}{*}{$\tau$} &\multicolumn{2}{c}{$\gamma=1$} &\multirow{2}{*}{$\tau$} &\multicolumn{2}{c}{$\gamma=1.5$} &\multirow{2}{*}{$\tau$}&\multicolumn{2}{c}{$\gamma=2$} \\
             \cline{3-4}          \cline{6-7}         \cline{9-10}
         &          &$e(N)$   &Order &         &$e(N)$   &Order &         &$e(N)$    &Order\\
\midrule
  40     &2.50e-02	&1.76e-01 &$-$   &4.77e-02	&4.12e-02 &$-$  &5.81e-02  &1.54e-02 &$-$\\
  80     &1.25e-02	&1.01e-01	&0.80  &2.48e-02	&1.79e-02	&1.27 &2.75e-02	 &6.20e-03	 &1.22\\
  160    &6.25e-03	&5.82e-02	&0.80  &1.24e-02	&7.80e-03	&1.19 &1.41e-02	 &2.64e-03	 &1.28\\
  320    &3.13e-03	&3.34e-02	&0.80  &6.33e-03	&3.40e-03	&1.24 &7.05e-03	 &1.14e-03	 &1.20\\
\midrule
\multicolumn{3}{l}{$\min\{\gamma\alpha, 2-\alpha\}$}   &0.80 & & &1.20 & & &1.20\\
\end{tabular*}
{\rule{\temptablewidth}{0.5pt}}	
\end{center}
\end{table}

\begin{table}[htb!]
\begin{center}
\caption{Temporal error of \eqref{eq: nonuniform L1 scheme} for $\alpha=0.4$ with $\gamma_{\mathrm{opt}}=4$.}\label{table:StabilizedL1-InitialSingularity-2}
\vspace*{0.3pt}
\def\temptablewidth{1.0\textwidth}
{\rule{\temptablewidth}{0.5pt}}
\begin{tabular*}{\temptablewidth}{@{\extracolsep{\fill}}cccccccccc}
\multirow{2}{*}{$N$} &\multirow{2}{*}{$\tau$} &\multicolumn{2}{c}{$\gamma=3$} &\multirow{2}{*}{$\tau$} &\multicolumn{2}{c}{$\gamma=4$} &\multirow{2}{*}{$\tau$}&\multicolumn{2}{c}{$\gamma=5$} \\
             \cline{3-4}          \cline{6-7}         \cline{9-10}
         &          &$e(N)$   &Order &         &$e(N)$   &Order &         &$e(N)$    &Order\\
\midrule
40   &5.95e-02    &5.03e-02    &$-$       &6.13e-02    &1.35e-02    &$-$     &6.65e-02    &9.27e-03    &$-$\\
80   &3.06e-02    &2.19e-02    &1.25      &3.02e-02    &4.44e-03    &1.56    &3.56e-02    &3.75e-03    &1.45\\
160  &1.55e-02    &9.54e-03    &1.22      &1.66e-02    &1.47e-03    &1.86    &1.64e-02    &1.18e-03    &1.49\\
320  &7.60e-03    &4.15e-03    &1.17      &8.49e-03    &4.88e-04    &1.64    &8.09e-03    &3.78e-04    &1.62\\
\midrule
\multicolumn{3}{l}{$\min\{\gamma\alpha, 2-\alpha\}$}   &1.20 & & &1.60 & & &1.60\\
\end{tabular*}
{\rule{\temptablewidth}{0.5pt}}
\end{center}
\end{table}	

\begin{example}\label{example:convergence test}
To calculate the errors in the mesh refinement tests,
we consider an exact solution $\Phi=\omega_{1+\alpha}(t)\sin(x)\sin(y)$ of
the TFMBE model with a proper forcing term $g(\mathbf{x}, t)$, i.e., $\partial_t^\alpha \Phi
+\kappa\brab{\epsilon^2\Delta^2 \Phi-\nabla\cdot f(\nabla \Phi)}=g(\mathbf{x}, t)$.
We solve it in the domain $\Omega=(0, 2\pi)^2$ with periodic boundary condition
by taking the model parameters $\kappa=1$ and $\epsilon=0.5$.
\end{example}

The spatial computational domain is divided into a $128^2$ uniform mesh.
The finial time is set as $T=1$.
We divided the time interval $[0, T]$ into two parts,
$[0, T_0]$ and $[T_0, T]$, with total $N$ subintervals.
In the interval $[0, T_0]$, we apply
the graded time mesh $t_k=(k/N_0)^\gamma$ for $0\le k\le N_0$,
where $T_0=\min\{1/\gamma, T\}$ and $N_0=\lceil\frac{N}{T+1-\gamma^{-1}}\rceil$.
The random time meshes with
$\tau_{N_0+k}:=(T-T_0)\varepsilon_k/S_1$ for $1\le k\le N_1$
are used in the remainder interval $[T_0, T]$
where $N_1:=N-N_0$, $S_1=\sum_{k=1}^{N_1} \varepsilon_k$ and
$\varepsilon_k \in(0, 1)$ are random numbers.

By setting different grading parameters $\gamma$, the numerical results in Table \ref{table:StabilizedL1-InitialSingularity-1}
and Table \ref{table:StabilizedL1-InitialSingularity-2}
are computed for the cases of $\alpha=0.8$  and $\alpha=0.4$, respectively.
It is seen from the tables that
when the graded parameters $\gamma<\gamma_{\mathrm{opt}}:=(2-\alpha)/\alpha$,
the L1 scheme \eqref{eq: nonuniform L1 scheme}
is of order $O(\tau^{\gamma \alpha})$. In addition, when $\gamma\ge \gamma_{\mathrm{opt}}$,
the optimal accuracy can reach to $O(\tau^{2-\alpha})$.
These results perfectly support the sharpness of our theoretical findings.

\subsection{Simulation of coarsening dynamics}
In this subsection, we will simulate the coarsening
dynamics of the TFMBE model.
We choose some appropriate adaptive time-stepping strategy and
depict the numerical behaviors of the original energy $E$
and the variational energy $\mathcal{E}_{\alpha}$ during the coarsening process.
\begin{example}\label{example:simulation of exact initial conditions}
We carry out a standard benchmark problem with the model parameters $\kappa=1$ and $\epsilon^2=0.1$,
and the initial data
$\phi(\mathbf{x}, 0)=0.1 \bra{\sin\bra{3x}\sin\bra{2y}+\sin\bra{5x}\sin\bra{5y}}.$
\end{example}

\begin{figure}[htb!]
\centering
\includegraphics[width=2in]{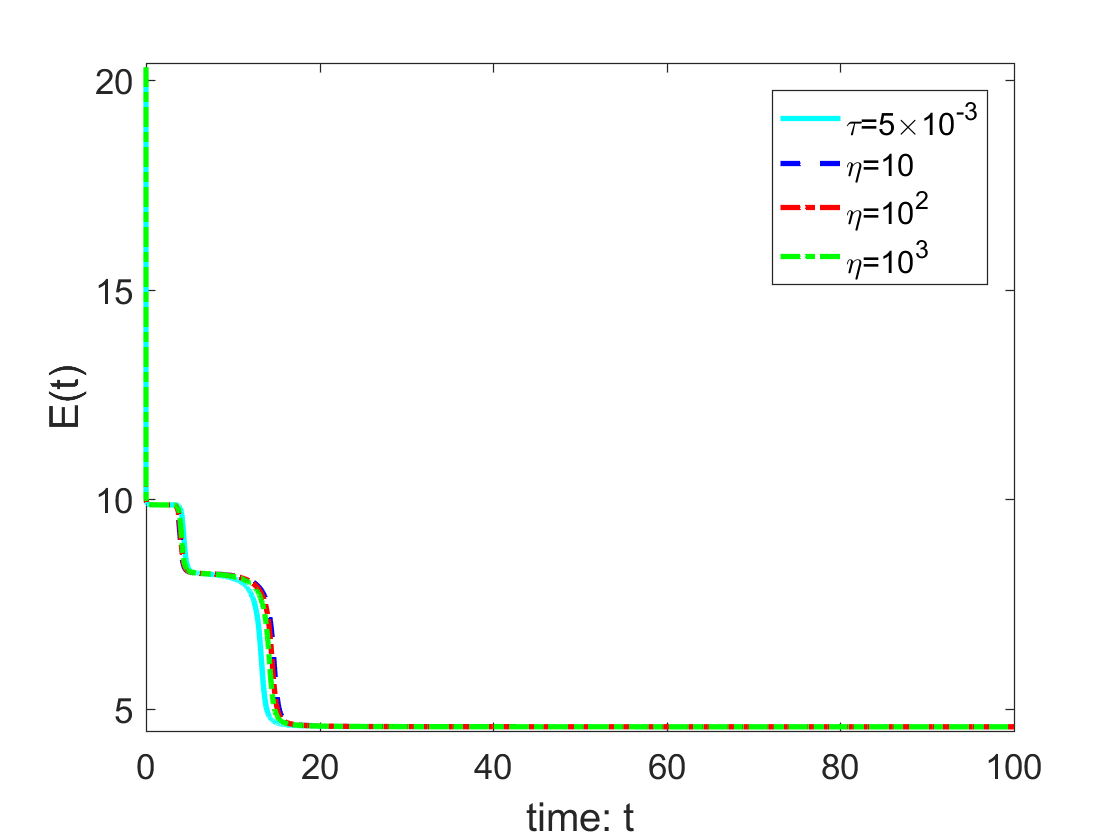}
\includegraphics[width=2in]{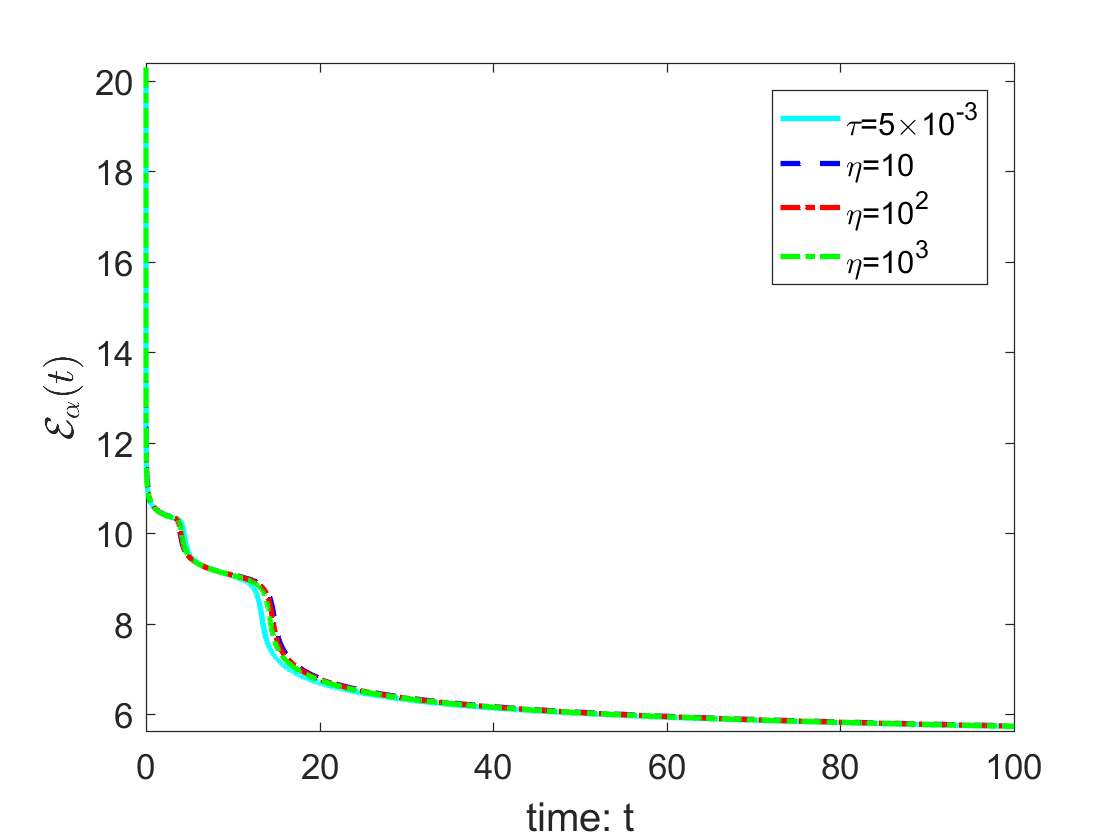}
\includegraphics[width=2in]{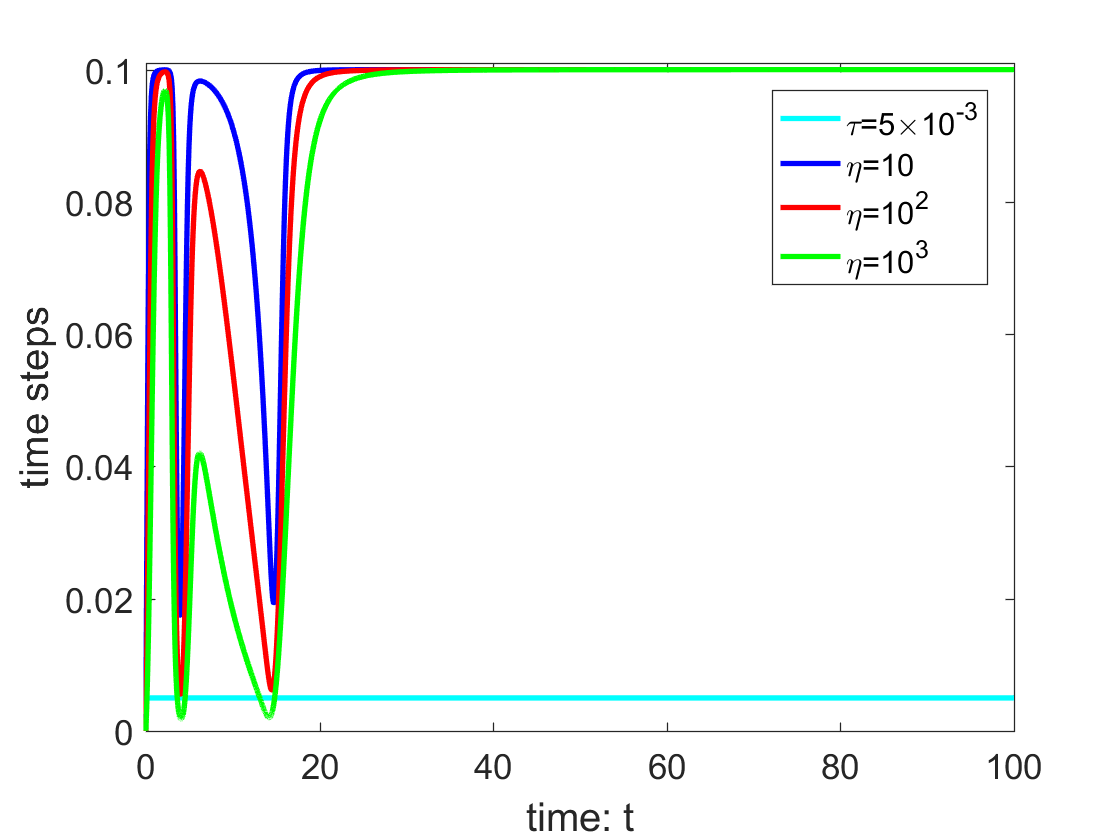}
\caption{Numerical comparisons of energy evolutions using uniform time step and adaptive
time-stepping strategy with different adaptive parameters $\eta$.}
\label{figure:etaEnergyEnalphaTimesteps}
\end{figure}

\begin{table}[htb!]
\begin{center}
\caption{Comparisons of CPU time (in seconds) and total time steps.}
\vspace*{0.3pt}
\def\temptablewidth{0.8\textwidth}
{\rule{\temptablewidth}{0.6pt}}
\begin{tabular*}{\temptablewidth}{@{\extracolsep{\fill}}c|cccccc}
Adaptive parameter    &$\eta=10$   &$\eta=10^2$   &$\eta=10^3$   &uniform step\\
\midrule
CPU time         &178.07        &222.81      &329.54     &1666.83\\
Time steps       &1156          &1496        &2669       &20048\\
\end{tabular*}
{\rule{\temptablewidth}{0.5pt}}
\label{table:CPU time}\vspace*{0.3pt}
\end{center}
\end{table}	

The TFMBE model has obvious multi-scale behaviors in time \cite{JiLiaoGongZhang:2018Adaptive}
and the variable-step L1 scheme \eqref{eq: nonuniform L1 scheme} is shown to be robustly stable and convergent
on arbitrary time meshes, see Theorems \ref{Thm:Energystability} and \ref{theorem:sharp L2 error estimate}.
So certain adaptive time-stepping approach is reasonably adopted in our numerical simulations
because it not only can capture the rapid changes of energy and numerical solution in a short time,
but also can improve the calculation efficiency with large time-steps when the solution varies slowly.

We select the time steps according to the change rate of the numerical solution with the following
adaptive time-stepping strategy, cf. \cite{LiaoZhuWang:2021},
\begin{align*}
\tau_{ada}=\max\bigg\{\tau_{\min}, \frac{\tau_{\max}}{\sqrt{1+\eta\mynormt{\partial_\tau \phi^n}^2}}\bigg\},
\end{align*}
where $\tau_{\max}$ and $\tau_{\min}$ are the predetermined maximum and minimum size of time-steps,
and $\eta$ is a user parameter to be determined.
The space domain $(0, 2\pi)^2$ is discretized by $128\times 128$ meshes during calculation.
In additional, let $\tau_{N_0}:=\tau_{\min}$
when the graded mesh is applied in the initial cell $[0, T_0]$
and the adaptive time-stepping strategy is employed in the remainder interval $[T_0, T]$,
in which $N_0$ is determined by $\tau_{N_0}=t_{N_0}-t_{N_0-1}$.

In order to determine a suitable parameter $\eta$,
we take $\tau_{\max}=10^{-1}$, $\tau_{\min}=10^{-3}$,
and consider  three different parameters
$\eta=10, 100$ and $10^3$.
The reference solution is computed by using the uniform time step $\tau=5\times 10^{-3}$.
As seen in Figure \ref{figure:etaEnergyEnalphaTimesteps},
the value of parameter $\eta$ evidently influences on the adaptive sizes of time steps.
Specially, when $\eta=10^3$, the time-steps have the smallest fluctuation,
and the L1 scheme can accurately capture the changes of original energy $E$ and
modified energy $\mathcal{E}_{\alpha}$ over the time.

\begin{figure}[htb!]
\centering
\includegraphics[width=1.47in]{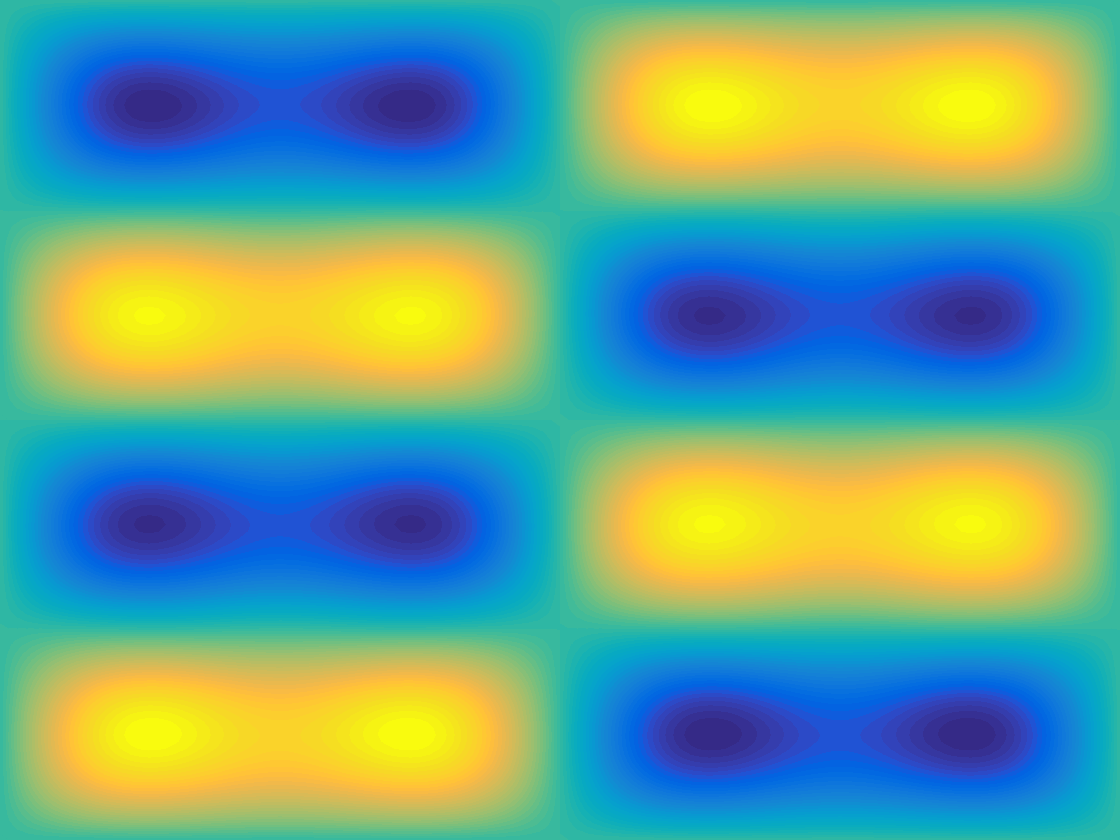}
\includegraphics[width=1.47in]{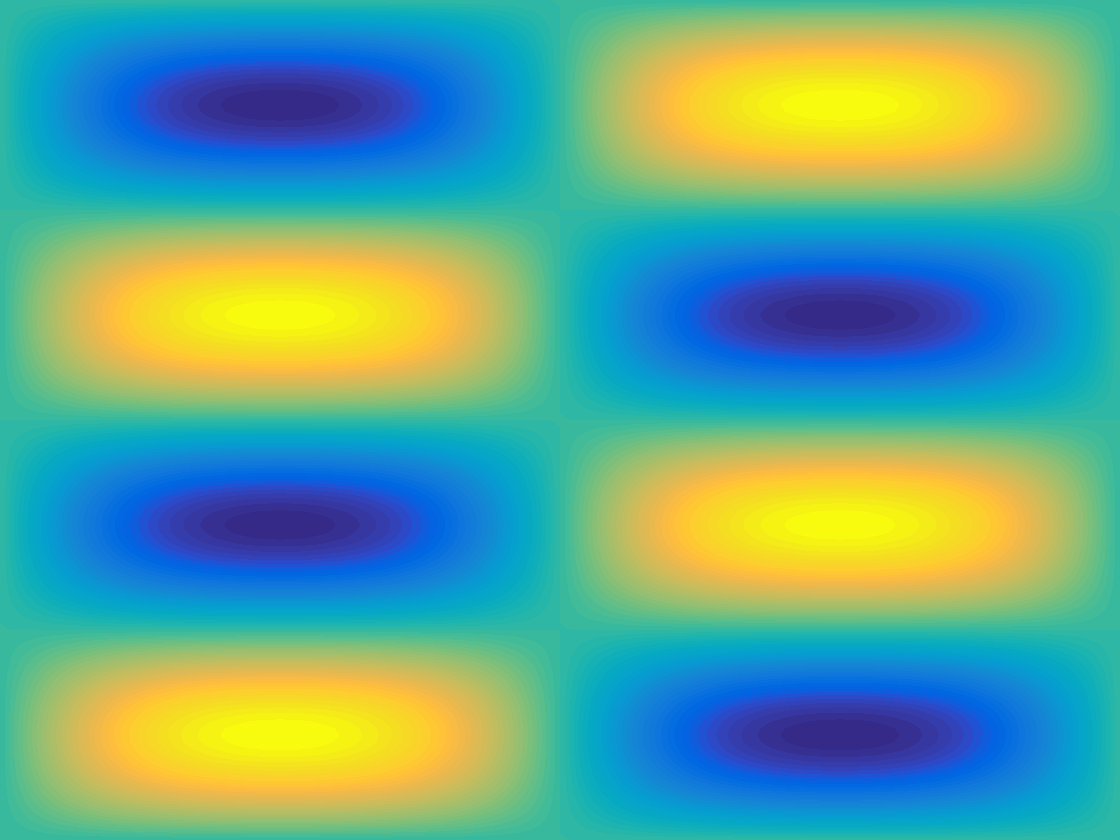}
\includegraphics[width=1.47in]{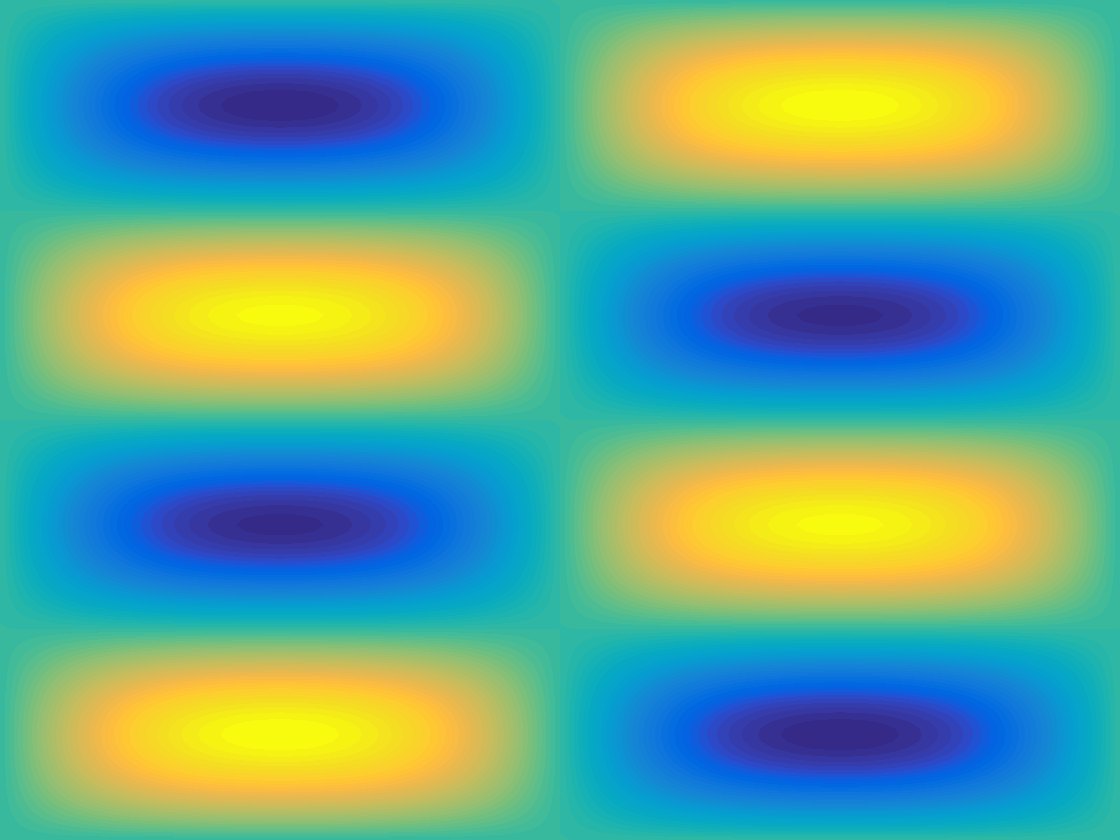}
\includegraphics[width=1.47in]{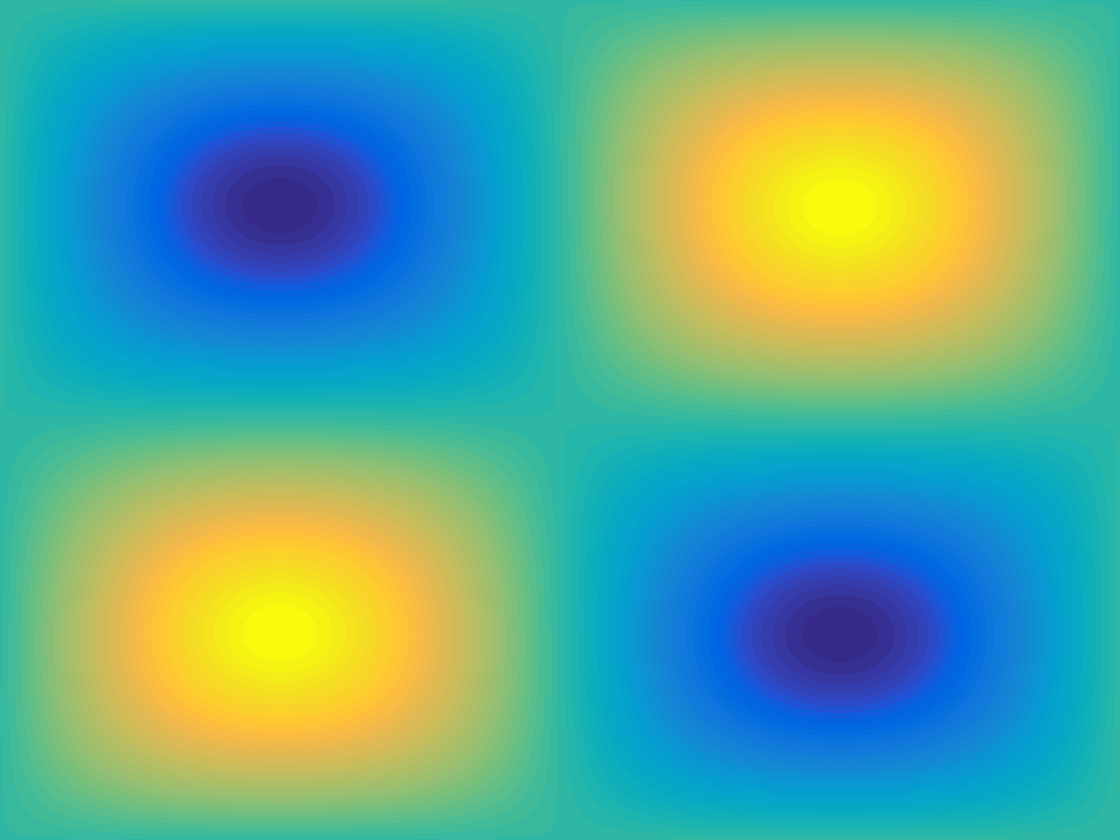}\\
\includegraphics[width=1.47in]{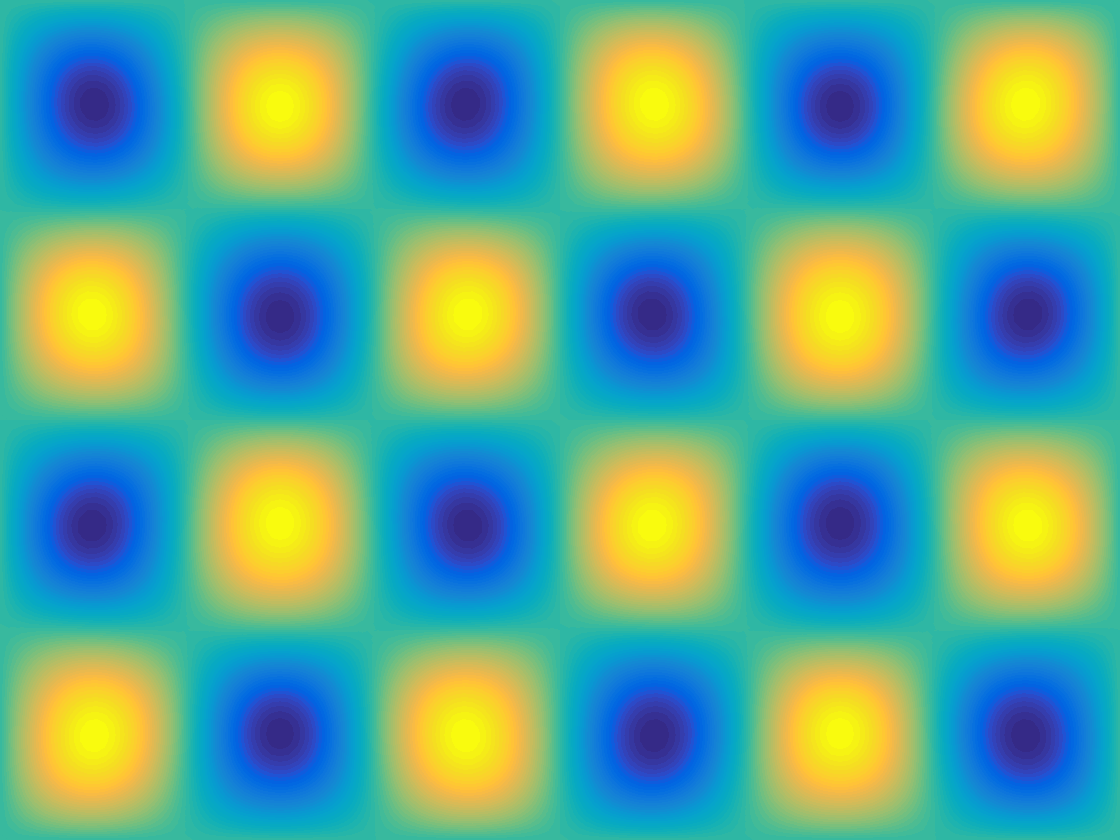}
\includegraphics[width=1.47in]{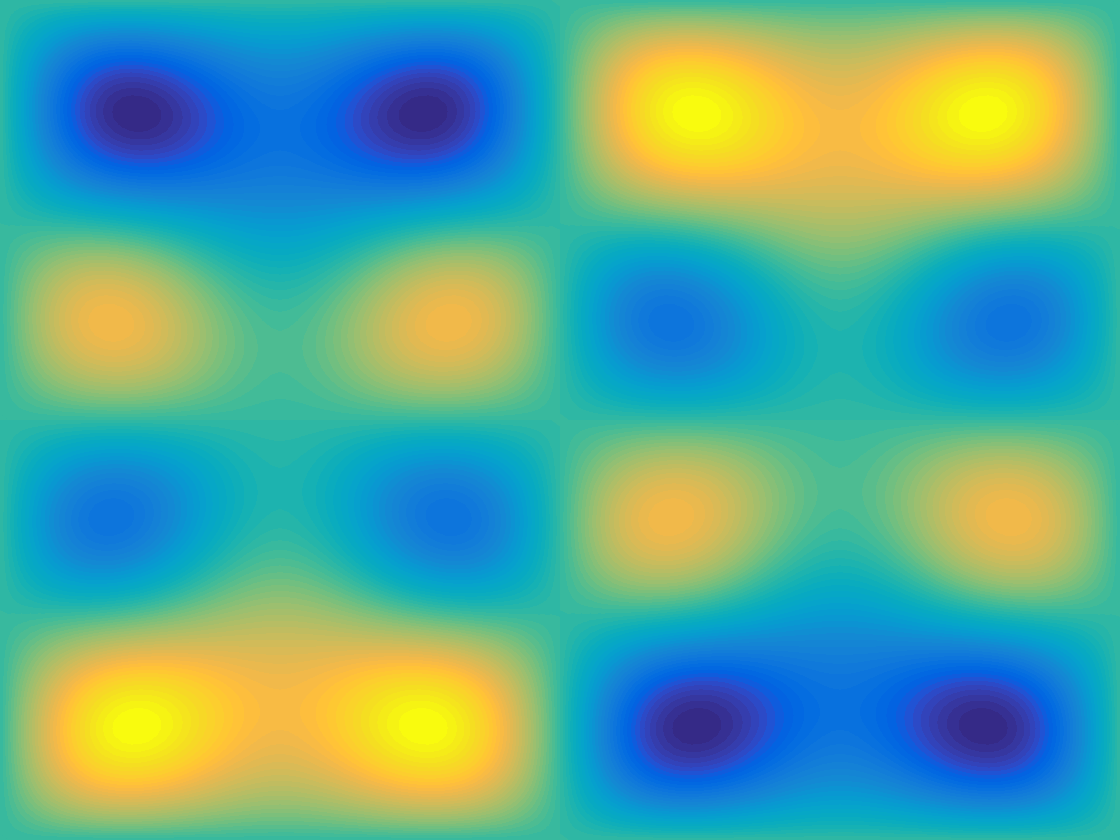}
\includegraphics[width=1.47in]{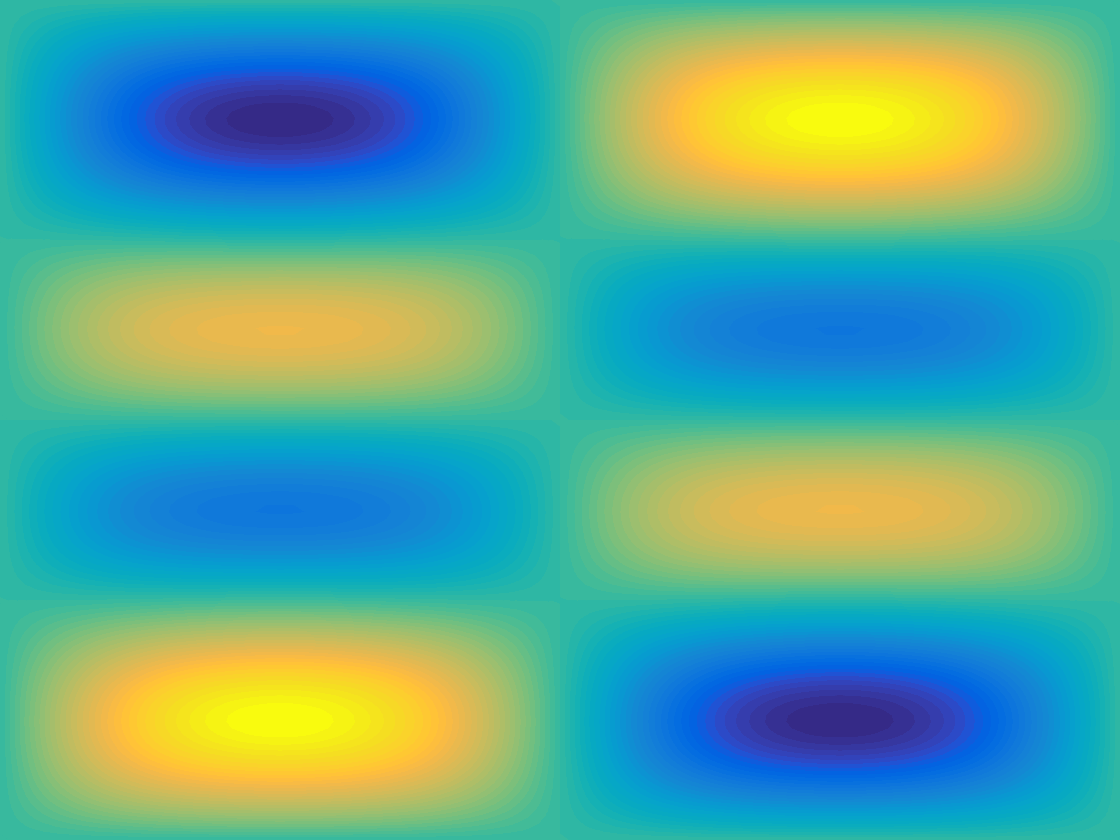}
\includegraphics[width=1.47in]{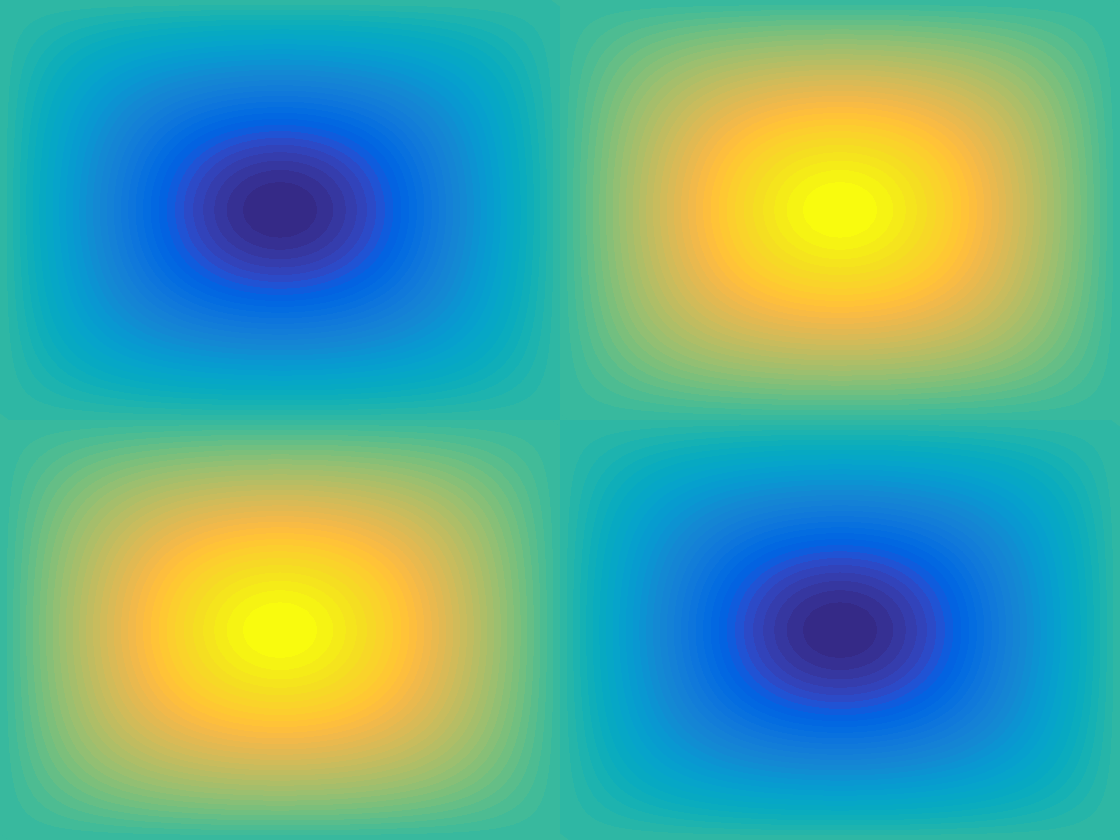}\\
\includegraphics[width=1.47in]{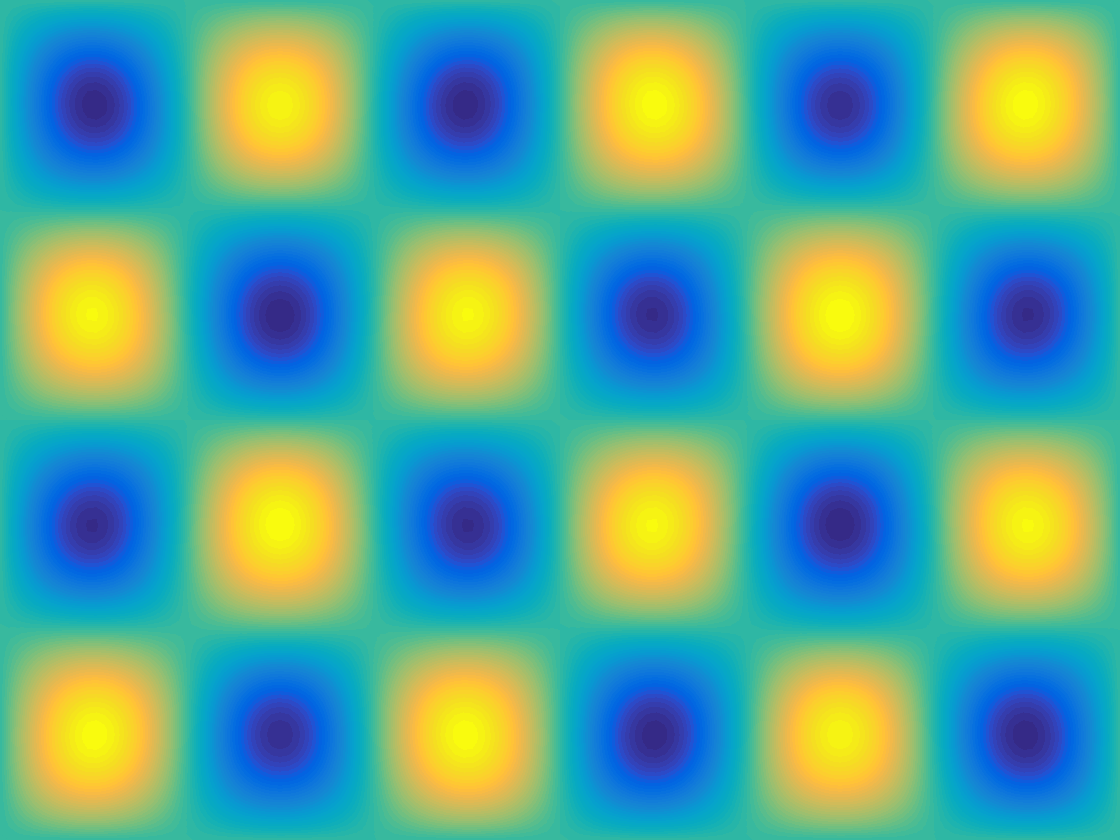}
\includegraphics[width=1.47in]{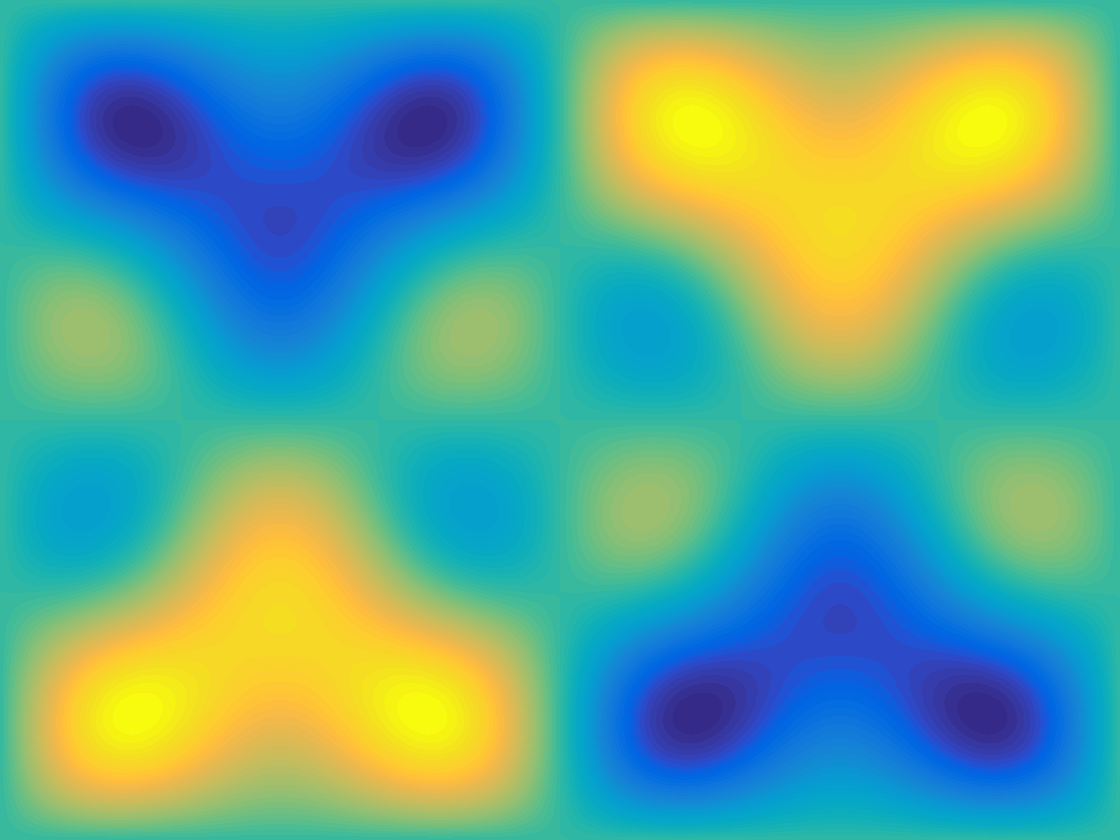}
\includegraphics[width=1.47in]{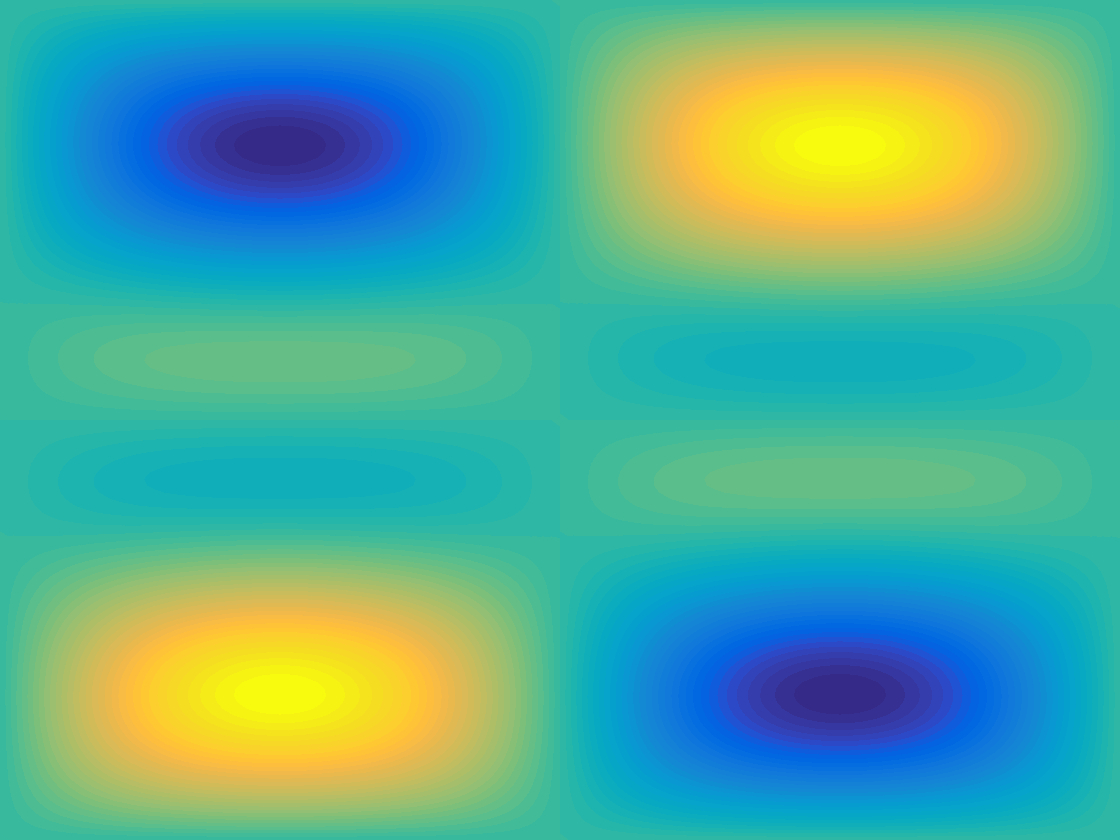}
\includegraphics[width=1.47in]{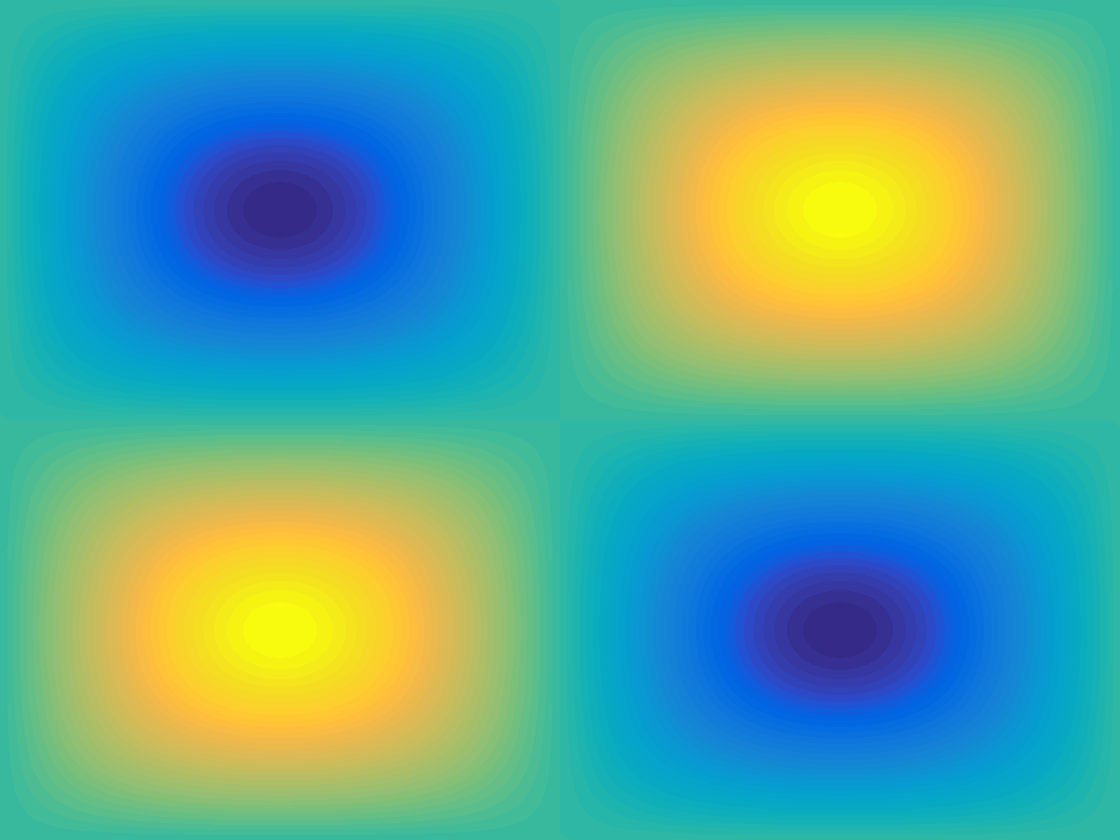}\\
\caption{Time snapshots of TFMBE model
\eqref{example:simulation of exact initial conditions} with $\varepsilon^2=0.1$ at
 $t=1.3, 3.0, 10, 50$ (from left to right)
for fractional orders $\alpha=0.4, 0.7, 0.9$ (from top to bottom), respectively.}
\label{figure:different alpha}
\end{figure}

\begin{figure}[htb!]
\centering
\includegraphics[width=2in]{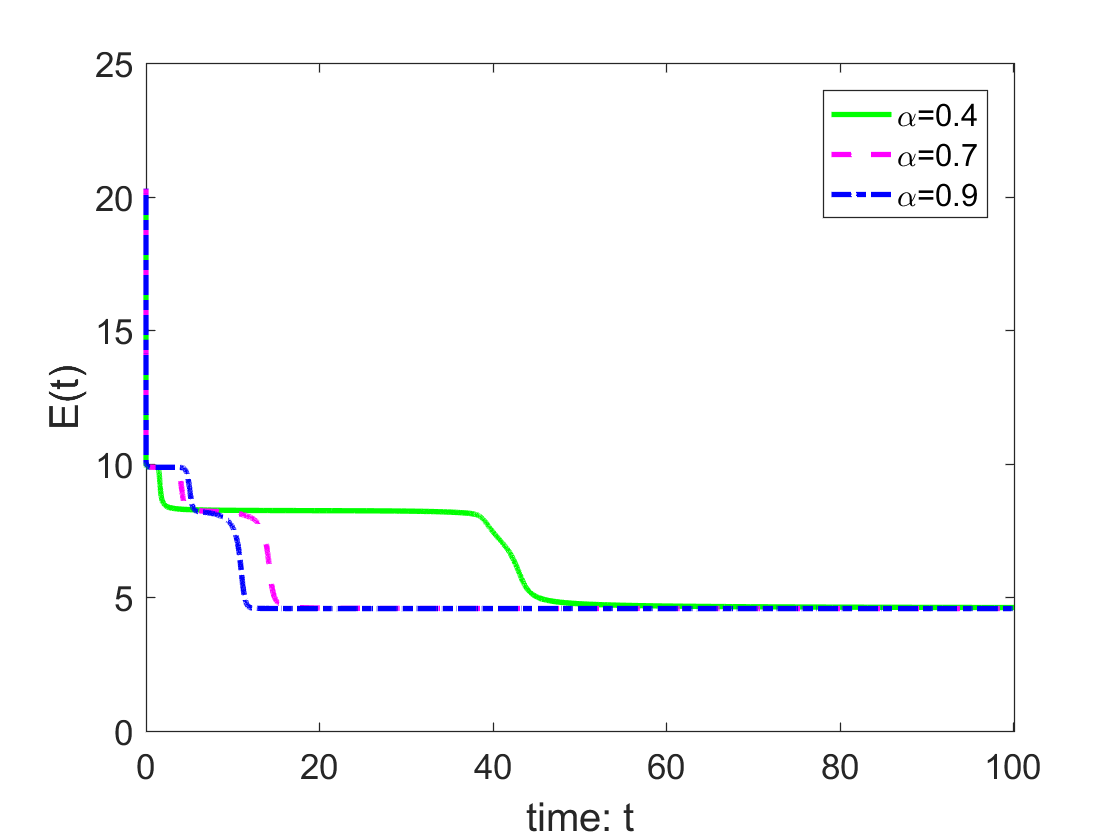}
\includegraphics[width=2in]{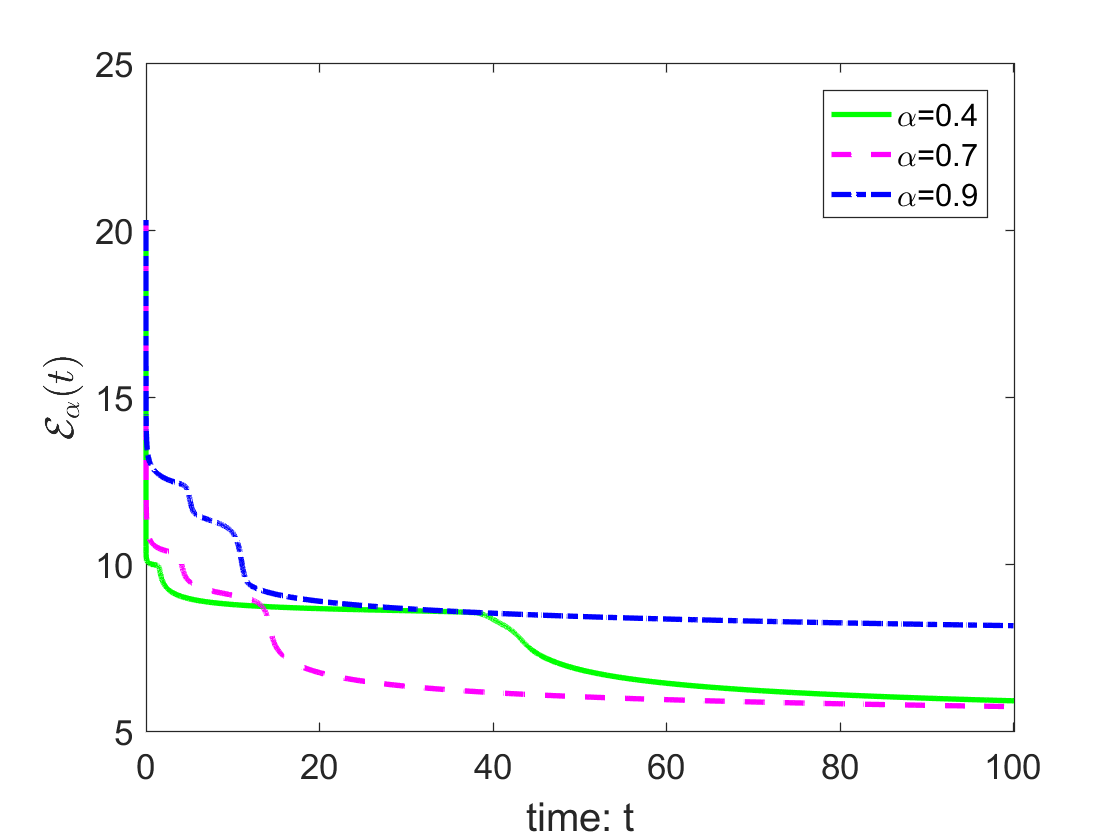}
\includegraphics[width=2in]{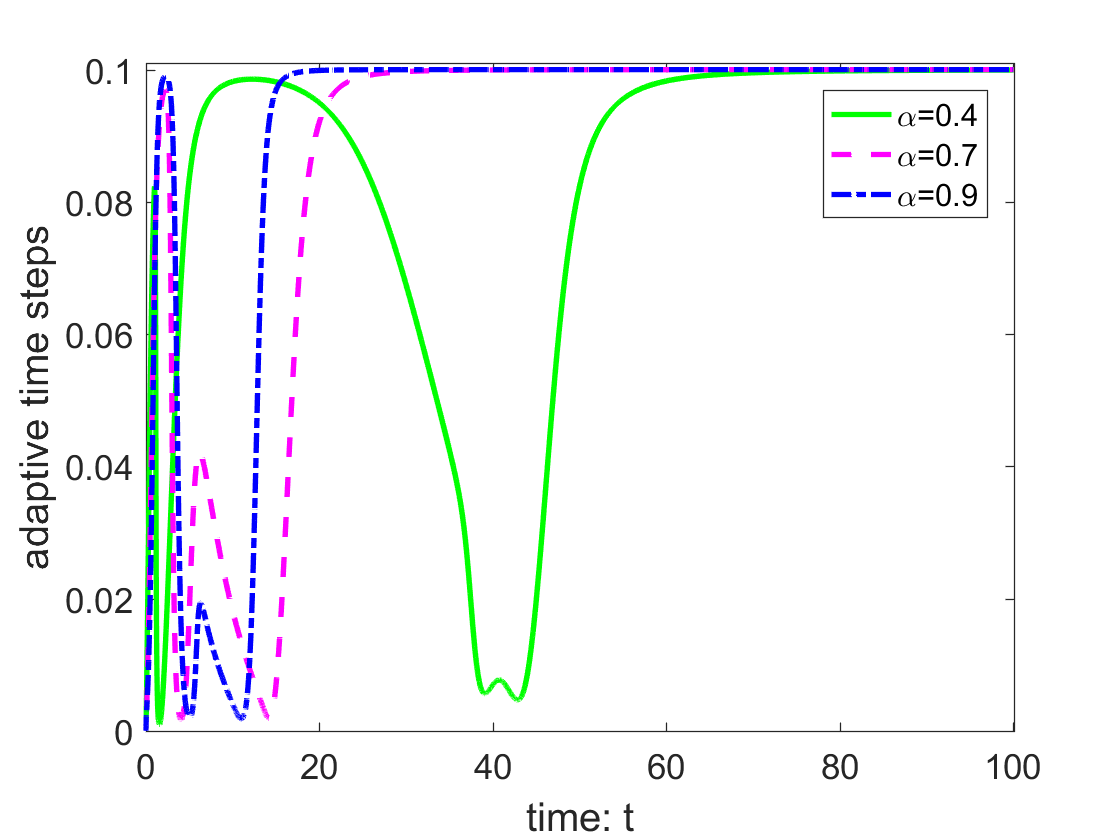}
\caption{Curves of original energy $E(t)$,
variational energy $\mathcal{E}_{\alpha}(t)$ and adaptive time steps $\tau_n$
generated for different fractional orders $\alpha$.}
\label{figures:adaptive}
\end{figure}

The corresponding CPU cost (in seconds) and the number of adaptive time levels
are listed in Table \ref{table:CPU time}.
We observe that, at least for this example, $\eta=10^3$ is a good choice
because it seems computationally more efficient than other cases using
the parameters $\eta=10$, $\eta=10^2$, and using the uniform step size.
As desired, the original energy $E$ monotonously decays over the time although we can not verify it theoretically.
On the other hand, as expected by our analysis,
the modified energy $\mathcal{E}_{\alpha}$ monotonously decays in the coarsening dynamics.

Next, by taking $\tau_{\max}=10^{-1}$, $\tau_{\min}=10^{-3}$ and
the parameter $\eta=10^3$ in the above adaptive time-stepping strategy,
we run the L1 scheme \eqref{eq: nonuniform L1 scheme}
for three different fractional orders $\alpha=0.4, 0.7$ and $0.9$ until the final time $T=100$.
The profiles of coarsening dynamics with different
fractional orders $\alpha=0.4,0.7$ and 0.9 for the TFMBE model \eqref{eq:TFMBE} are shown
in Figure \ref{figure:different alpha},
where the snapshots of solution profiles are taken at time $t=1.3, 3.0, 10$ and 50, respectively.
We observe that the coarsening rates are always dependent on the fractional order and the time period,
but they all approach the steady state near $t=50$.
The curves of original energy $E$ and the variational energy $\mathcal{E}_{\alpha}$
over the time interval $t\in[0,100]$
are depicted in Figure
\ref{figures:adaptive}. The initial energy decays rapidly in all cases,
while it decays slower for the smaller fractional order $\alpha$.
As the time goes on, the evolution dynamics reach
the same steady state in the end for different fractional orders.
These results are in accordance with the previous observations in
\cite{ChenZhaoCaoWangZhang:2018,JiLiaoGongZhang:2018Adaptive}.


\section{Acknowledgements}
The authors would like to thank Prof. Gong Yuezheng, Dr. Ji Bingquan and Ms. Zhu Xiaohan for their valuable suggestions.

\end{document}